\documentclass[12pt]{article}

\usepackage{amsmath,amssymb,amsthm,mathtools,bbm}
\usepackage{amsfonts}
\usepackage{bigints}
\usepackage{enumerate}
\usepackage{hyperref}
\usepackage{mathtools}

\DeclareMathOperator{\var}{var}

\newcommand{\R}{\mathbb{R}}
\newcommand{\N}{\mathbb{N}}

\newcommand{\dimh}{\dim_{\mathrm H}}
\newcommand{\dime}{\dim_{\mathrm e}}
\newcommand{\Dn}{\mathcal D_n}
\newcommand{\Dl}{\mathcal D_\ell}

\newcommand{\cF}{\mathcal F}
\newcommand{\cP}{\mathcal P}
\newcommand{\cB}{\mathcal B}
\newcommand{\cQ}{\mathcal Q}

\newcommand{\tbu}{\textbf{u}}
\newcommand{\tbv}{\textbf{v}}

\begin{document}

\newcommand{\A}{\mbox{${{{\cal A}}}$}}

%1111111111111111111111111111111111111111111111111111111111111111111111111111

\author{Attila Losonczi}
\title{Lower semicontinuity of dimension for ergodic and quasi-Bernoulli measures under conformal projections}

\date{\today}%\date{22 September 2026}

\newtheorem{thm}{\qquad Theorem}[section]
\newtheorem{prp}[thm]{\qquad Proposition}
\newtheorem{lem}[thm]{\qquad Lemma}
\newtheorem{cor}[thm]{\qquad Corollary}
\newtheorem{rem}[thm]{\qquad Remark}
\newtheorem{ex}[thm]{\qquad Example}
\newtheorem{prb}[thm]{\qquad Problem}
\newtheorem{df}[thm]{\qquad Definition}
\newtheorem{fact}[thm]{Fact}

\maketitle

\begin{abstract}
We study the Hausdorff dimension of projections of symbolic measures under the coding maps of conformal iterated function systems on the real line, without separation assumptions. For each fixed shift-invariant ergodic measure, its projection entropy and the Hausdorff dimension of its projection are lower semicontinuous under $C^1$ perturbations of the system. When the symbolic measure also varies, we give a sharp bound on the possible loss of projected dimension in terms of the loss of symbolic entropy and the Lyapunov exponent of the limiting system. This yields joint lower semicontinuity for weakly convergent quasi-Bernoulli measures with a common quasi-Bernoulli constant, including measures that are not shift-invariant. For invariant Gibbs measures associated with Hölder potentials, uniform convergence of the potentials suffices. Finally, for each fixed shift-invariant quasi-Bernoulli measure, we obtain the entropy-to-Lyapunov dimension formula, truncated at one, on a dense $G_\delta$ set of systems whose derivative moduli are strictly less than $1/2$.

\noindent
\footnotetext{\noindent
MSC 2020: Primary 28A80; Secondary 37C45, 37D35 

\noindent Key Words and Phrases: Hausdorff dimension; entropy dimension; lower semicontinuity; conformal iterated function systems; quasi-Bernoulli measures; Gibbs measures.
}

\end{abstract}

%\tableofcontents

%-------------------------------------------------------------------------------------------------------------------------------------------Introduction----------------------------------
\section{Introduction}\label{sec:intro}
The dimension theory of iterated function systems studies how the
geometry of the defining maps and the statistical properties of
symbolic measures determine the size of their projections. For
conformal systems, entropy and Lyapunov exponents provide a natural
upper bound for the Hausdorff dimension of a projected invariant
measure. Overlaps between cylinder sets may prevent this bound from
being attained, and they also make the dependence of dimension on the
defining system more difficult to describe. A basic stability question
is whether the dimension of a projected measure is lower
semicontinuous when the maps are perturbed. We address this question
for conformal iterated function systems on the real line and fixed
quasi-Bernoulli measures on the symbolic space.

Entropy is an important tool for studying dimension in the presence
of overlaps. Peres and Solomyak \cite{PeresSolomyak00} proved the
existence of \(L^q\)-dimensions and entropy dimension for self-conformal
measures without separation assumptions. Their proof uses coverings
adapted to the geometry of cylinder sets and comparisons between
entropy at different scales. Ngai and Xu \cite{NgaiXu23} obtained
related existence results on Riemannian manifolds; their
entropy-dimension theorem assumes that the self-conformal measure
is doubling. A complementary approach was developed by Feng and
Hu \cite{FengHu09}, who introduced projection entropy and proved
exact dimensionality for conformal projections of invariant ergodic
measures, again without separation assumptions. These results form
part of the broader connection between dimension, entropy, and
Lyapunov exponents studied by Young \cite{Young82}.

%For self-similar Bernoulli measures, lower semicontinuity of dimension is treated by the entropy method in\cite[Section~3.7]{BSS23}. 
For self-similar Bernoulli measures, lower semicontinuity of
Hausdorff dimension follows from the entropy estimates of
Peres and Solomyak \cite{PeresSolomyak00}, together with the
exact-dimensionality theorem of Feng and Hu \cite{FengHu09}.
An exposition of the entropy argument is given in
\cite[Section~3.7]{BSS23}. Extending this approach to nonlinear
conformal maps requires distortion estimates that are uniform over
all cylinder lengths. Passing from Bernoulli to quasi-Bernoulli
measures introduces a further issue: after conditioning on a symbolic
cylinder and removing its prefix, the resulting measure need not
coincide with the original symbolic measure. The quasi-Bernoulli
property supplies a uniform comparison between these conditional
measures and the original measure, and this comparison replaces the
exact product structure in our argument.

Quasi-Bernoulli measures include Bernoulli measures and Gibbs measures
for H\"older continuous potentials on a finite full shift. They
therefore allow dependence between successive symbols while retaining
uniform multiplicative control of cylinder probabilities. We refer
to Heurteaux \cite{Heurteaux2007} for their dimension-theoretic and
ergodic properties, and to Walters \cite{Walters82} for the underlying
ergodic theory and thermodynamic formalism. Related multiplicative
estimates also occur in the almost-additive thermodynamic formalism
for matrix cocycles studied by B\'ar\'any, K\"aenm\"aki, and
Morris \cite{BKM2020}.

We study how the Hausdorff dimension of a projected measure changes
when the defining system and the symbolic measure are perturbed.
Let $\Pi_{\mathcal F}$ denote the coding map of the system $\mathcal F$.
For a fixed shift-invariant ergodic measure $\mu$, we prove that the
projection entropy $h_{\Pi_{\mathcal F}}(\sigma,\mu)$ is lower
semicontinuous as $\mathcal F$ varies in the $C^1$ topology
(Theorem~\ref{thm:lsc-ergodic}). The proof expresses projection entropy
in terms of the conditional entropies of finite cylinder partitions
given the coding map. Since the Lyapunov exponent depends continuously
on $\mathcal F$, the dimension formula of Feng and Hu then gives lower
semicontinuity of $\dim_H((\Pi_{\mathcal F})_*\mu)$.

When the symbolic measure also varies, weak convergence alone allows
the projected dimension to drop, even for a fixed pair of similarities.
Our quantitative result relates this drop to the loss of symbolic
entropy. Suppose that $\mathcal F_i\to\mathcal F$ in $C^1$ and that
$\mu_i\Rightarrow\mu$, where all the measures are shift-invariant and
ergodic. Writing $d_i=\dim_H((\Pi_{\mathcal F_i})_*\mu_i)$ and
$d=\dim_H((\Pi_{\mathcal F})_*\mu)$, we prove that
\[
  d-\liminf_{i\to\infty}d_i
  \leq
  \frac{h_\mu(\sigma)-\liminf_{i\to\infty}h_{\mu_i}(\sigma)}
       {\chi_\mu(\mathcal F)}
  \qquad\text{(Theorem~\ref{thm:lsc-joint}).}
\]
In particular, convergence of the symbolic entropies implies joint
lower semicontinuity of the projected dimension. The bound is sharp:
for the binary expansion maps, periodic orbit measures can converge
weakly to the fair Bernoulli measure, with projected dimension zero
along the sequence and one at the limit.

For quasi-Bernoulli measures, we give a direct proof based on
finite-scale entropy estimates. If $\mu_i\Rightarrow\mu$ with a common
quasi-Bernoulli constant and $\mathcal F_i\to\mathcal F$ in
$C^{1+\theta}$, then the projected dimension is jointly lower
semicontinuous (Theorem~\ref{thm:gen}). Here the symbolic measures
need not be shift-invariant. For invariant Gibbs measures associated
with H\"older continuous potentials, the entropy bound above gives
a further conclusion: $C^1$ convergence of the systems and uniform
convergence of the potentials suffice for joint lower semicontinuity.
Indeed, the variational principle yields convergence of the symbolic
entropies in this setting. None of these semicontinuity results
requires a separation condition.

As an application, we obtain a generic dimension formula for systems
whose derivative moduli are strictly less than $1/2$. For each fixed
shift-invariant quasi-Bernoulli measure $\mu$, there is a dense
$G_\delta$ set of such systems for which the projected dimension equals
$\min\{1,h_\mu(\sigma)/\chi_\mu(\mathcal F)\}$
(Theorem~\ref{thm:residual}). The set may depend on $\mu$.
Here lower semicontinuity is combined with transversality for
translation families as presented in \cite[Chapter~14]{BSS23}.

The remainder of this section introduces the notation and the
classes of systems and measures under consideration. Section~2
recalls the covering and entropy estimates used in the proofs,
establishes lower semicontinuity, and derives the generic
dimension formula.

%------------------------------------------------------------------------------------------------------------------------------Basic notations---------------------------
\subsection{Preliminaries and notation}

We follow the general notation and terminology of the book of Bárány, Simon and Solomyak \cite{BSS23}.
All measures are Borel probability measures unless otherwise stated.
All logarithms are natural, and \(0\log0=0\). For a bounded
nonempty set \(A\subset\mathbb R\), we write
\(|A|=\operatorname{diam}A\); in particular, \(|I|\) is the length
of an interval \(I\). We use \(\mathbb N=\{1,2,\ldots\}\).

\paragraph{Conformal iterated function systems.}
Fix a nondegenerate compact interval \(I\subset\mathbb R\), an
integer \(m\geq2\), and \(\theta\in(0,1)\). Let \(\Theta(I)\)
be the collection of ordered systems
\(\mathcal F=(f_1,\ldots,f_m)\) satisfying
\[
 f_i\in C^{1+\theta}(I),\qquad
 f_i(I)\subset\operatorname{int}(I),\qquad
 0<|f_i'(x)|<1
 \quad(x\in I,\ 1\leq i\leq m).
\]
Both orientation-preserving and orientation-reversing maps are
allowed. For \(0<\gamma_1<\gamma_2<1\), set
\[
 \Theta_{\gamma_1,\gamma_2}(I)
 =\left\{\mathcal F\in\Theta(I):
 \gamma_1\leq|f_i'(x)|\leq\gamma_2
 \text{ for all }x\in I,\ 1\leq i\leq m\right\}.
\]
We also use the separate notation
\[
 \Theta_{0,1/2}(I)
 =\left\{\mathcal F\in\Theta(I):
 |f_i'(x)|<\tfrac12
 \text{ for all }x\in I,\ 1\leq i\leq m\right\}.
\]
The strict inequality in the last definition is intentional.

For \(g\in C^\theta(I)\), write
\[
 [g]_\theta=\sup_{\substack{x,y\in I\\x\ne y}}
 \frac{|g(x)-g(y)|}{|x-y|^\theta}.
\]
We equip \(C^{1+\theta}(I)\) with the norm
\[
 \|f\|_{1+\theta}
 =\|f\|_\infty+\|f'\|_\infty+[f']_\theta,
\]
and \(\Theta(I)\), together with its subspaces, with the metric
\[
 \rho(\mathcal F,\mathcal G)
 =\max_{1\leq i\leq m}\|f_i-g_i\|_{1+\theta}.
\]
We write \(L(\mathcal F)=\max_i[f_i']_\theta\).

\paragraph{Symbolic space and coding.}
Let \(\Sigma=\{1,\ldots,m\}^{\mathbb N}$, $\Sigma_n=\{1,\ldots,m\}^n\), and
\(\Sigma_*=\bigcup_{k=0}^\infty\Sigma_k\), where
\(\Sigma_0=\{\emptyset\}\). The left shift is denoted by
\(\sigma\). If \(\mathbf u=u_1\cdots u_k\), then
\(|\mathbf u|=k\), \(\mathbf u^-=u_1\cdots u_{k-1}\) for
\(k\geq1\), and
\[
 [\mathbf u]=\{\omega\in\Sigma:\omega_1=u_1,\ldots,\omega_k=u_k\}.
\]
We set \([\emptyset]=\Sigma\), and write \(\mathbf u\mathbf v\)
for concatenation and \(\omega|_k=\omega_1\cdots\omega_k\) for $\omega\in\Sigma$.

For \(\mathcal F\in\Theta(I)\), define
\[
 f_{\mathbf u}=f_{u_1}\circ\cdots\circ f_{u_k},
 \qquad f_\emptyset=\mathrm{id}_I.
\]
The coding map is
\[
 \Pi_{\mathcal F}(\omega)
 =\lim_{k\to\infty}f_{\omega|_k}(x_0),
 \qquad \omega\in\Sigma,
\]
where \(x_0\in I\) is arbitrary. The limit is independent of
\(x_0\). Its image is the attractor
\(\Lambda_{\mathcal F}=\Pi_{\mathcal F}(\Sigma)\), the unique
nonempty compact set satisfying
\[
 \Lambda_{\mathcal F}=\bigcup_{i=1}^m f_i(\Lambda_{\mathcal F}).
\]
We write
\[
 I_{\mathbf u}=f_{\mathbf u}(I),\qquad
 \Lambda_{\mathbf u}=f_{\mathbf u}(\Lambda_{\mathcal F}),
 \qquad \Lambda_\emptyset=\Lambda_{\mathcal F}.
\]
The dependence on \(\mathcal F\) is omitted when it is clear.

The bounded distortion property gives a constant
\(D(\mathcal F)\geq1\) such that
\[
 D(\mathcal F)^{-1}
 \leq\frac{|f_{\mathbf u}'(x)|}{|f_{\mathbf u}'(y)|}
 \leq D(\mathcal F)
 \quad(\mathbf u\in\Sigma_*,\ x,y\in I).
\]
On \(\Theta_{\gamma_1,\gamma_2}(I)\), one may take
\[
 D(\mathcal F)=
 \exp\left(
 \frac{L(\mathcal F)|I|^\theta}
 {\gamma_1(1-\gamma_2^\theta)}
 \right).
\]
In particular, distortion constants can be chosen uniformly in a
sufficiently small \(C^{1+\theta}\) neighbourhood of a fixed system.

\paragraph{Quasi-Bernoulli measures and their projections.}
A probability measure \(\mu\) on \(\Sigma\) is
\emph{quasi-Bernoulli} if there exists \(K\geq1\) such that
\[
 K^{-1}\mu([\mathbf u])\mu([\mathbf v])
 \leq\mu([\mathbf u\mathbf v])
 \leq K\mu([\mathbf u])\mu([\mathbf v])
 \quad(\mathbf u,\mathbf v\in\Sigma_*).
\]
Shift-invariance is not part of this definition.
For a probability vector \(p=(p_1,\ldots,p_m)\), the Bernoulli
measure \(\mu_p=p^{\mathbb N}\) satisfies this condition with
\(K=1\).

The projection of \(\mu\) under \(\mathcal F\) is
\[
 \nu_{\mathcal F}=(\Pi_{\mathcal F})_*\mu,
 \qquad
 \nu_{\mathcal F}(A)=\mu(\Pi_{\mathcal F}^{-1}(A)).
\]
%When \(\mu=\mu_p\), we also write \(\nu_{\mathcal F,p}\);
%this is the self-conformal measure satisfying
%\[
% \nu_{\mathcal F,p}=\sum_{i=1}^m p_i(f_i)_*\nu_{\mathcal F,p}.
%\]
For a general quasi-Bernoulli measure, we refer to
\(\nu_{\mathcal F}\) as a conformal projection of \(\mu\).

%Every quasi-Bernoulli measure \(\mu\) has a shift-invariant,
%ergodic, quasi-Bernoulli representative \(\widetilde\mu\)
%with
%\[
% C^{-1}\mu\leq\widetilde\mu\leq C\mu
%\]
%for some \(C\geq1\); see \cite[Section~5]{Heurteaux2007}.
%If \(\mu\) is invariant, we take \(\widetilde\mu=\mu\).
%The invariant representative is unique among invariant
%probability measures equivalent to \(\mu\), since equivalent
%invariant ergodic probability measures coincide.
%For every \(\mathcal F\), the corresponding projections satisfy
%\[
% C^{-1}\nu_{\mathcal F}
% \leq(\Pi_{\mathcal F})_*\widetilde\mu
% \leq C\nu_{\mathcal F}.
%\]

%\paragraph{Entropy and Lyapunov exponent.}
%For a quasi-Bernoulli measure \(\mu\), define
%\[
% h_\mu=h_{\widetilde\mu}(\sigma)
% =\lim_{k\to\infty}
% -\frac1k\sum_{\mathbf u\in\mathcal A^k}
% \widetilde\mu([\mathbf u])
% \log\widetilde\mu([\mathbf u]),
%\]
%and
%\[
% \chi_\mu(\mathcal F)
% =-\int_\Sigma
% \log\left|f_{\omega_1}'\bigl(\Pi_{\mathcal F}(\sigma\omega)\bigr)\right|
% \,d\widetilde\mu(\omega).
%\]
%Thus the notation agrees with the usual entropy and Lyapunov
%exponent when \(\mu\) is invariant. For \(\mu=\mu_p\), we
%have \(h_{\mu_p}=-\sum_i p_i\log p_i\).
%On \(\Theta_{\gamma_1,\gamma_2}(I)\),
%\[
% -\log\gamma_2\leq\chi_\mu(\mathcal F)\leq-\log\gamma_1.
%\]

\paragraph{Dimensions and dyadic entropy.}
For a compactly supported probability measure \(\nu\) on
\(\mathbb R\), set
\[
 \dim_H\nu
 =\inf\{\dim_H A:A\subset\mathbb R\text{ is Borel},\ \nu(A)=1\}.
\]
The measure is \emph{exact dimensional} with dimension \(d\)
if
\[
 \lim_{r\downarrow0}\frac{\log\nu(B(x,r))}{\log r}=d
 \quad\text{for }\nu\text{-almost every }x.
\]
In this case \(\dim_H\nu=d\).

For \(n\in\mathbb Z\), let
\[
 \mathcal D_n
 =\{[j2^{-n},(j+1)2^{-n}):j\in\mathbb Z\},
\]
and denote by \(\mathcal D_n(x)\) the unique member containing
\(x\). Define
\begin{equation}\label{eq:dyadic-entropy}
 h_n(\nu)%=H(\nu,\mathcal D_n)
 =-\sum_{Q\in\mathcal D_n}\nu(Q)\log\nu(Q).
\end{equation}
The entropy dimension, when the limit exists, is
\[
 \dim_{\mathrm e}\nu
 =\lim_{n\to\infty}\frac{h_n(\nu)}{n\log2}.
\]
We will also use translated nested grids
\(\mathcal D_n^a=\{a+Q:Q\in\mathcal D_n\}\), with a single
translation \(a\) fixed for all \(n\). For any given \(\nu\),
one can choose \(a\) so that every grid boundary has zero
\(\nu\)-measure. Passing between two such grids changes each
entropy by at most \(\log2\), so it does not affect entropy
dimension.

By the exact-dimensionality theorem of Feng and
Hu \cite[Theorem~2.8]{FengHu09}, the projection of an invariant
ergodic measure under a conformal IFS is exact dimensional.
%For non-invariant quasi-Bernoulli measures, the same conclusion follows by comparison with \(\widetilde\mu\).
For these compactly supported projections, exact dimensionality
also gives
\[
 \dim_{\mathrm e}\nu_{\mathcal F}
 =\dim_H\nu_{\mathcal F};
\]
see \cite{Young82} and the dyadic local-dimension framework
of \cite{KaenmakiRajalaSuomala2010}.

\paragraph{Hölder continuous functions.}
Fix $0<r<1$ and define the metric on $\Sigma$
\[d_r(\tau,\omega)=r^{|\tau\wedge\omega|}\]
where $\tau\wedge\omega$ is their common initial word, and $|\tau\wedge\omega|$ is its length.

If $\varphi:\Sigma\to\R$, \(0<\alpha<1\) set
\[[\varphi]_\alpha=\sup_{\tau\ne\omega}\frac{|\varphi(\tau)-\varphi(\omega)|}{d_r(\tau,\omega)^\alpha},\]
which is the smallest possible Hölder constant, called the Hölder seminorm. The $C^\alpha$ norm is
\[||\varphi||_\alpha=||\varphi||_\infty+[\varphi]_\alpha.\]

If $\varphi:\Sigma\to\R$, set
\[\var_n(\varphi)=\sup\{|\varphi(\tau)-\varphi(\omega)|:|\tau\wedge\omega|=n\},\]  %=\sup\{|\varphi(\tau)-\varphi(\omega)|:\tau|_n=\omega|_n\}
\[D_\varphi=\sum_{n=1}^\infty\var_n(\varphi).\]

If $\varphi$ is $\alpha$-Hölder with Hölder seminorm $[\varphi]_\alpha$, then
\[\var_n(\varphi)\leq [\varphi]_\alpha r^{\alpha n},\]
hence
\[D_\varphi\leq [\varphi]_\alpha\frac{r^\alpha}{1-r^\alpha}<\infty.\]

If $\varphi_i\to\varphi$ in $C^\alpha(\Sigma)$ according to the previous norm, then
\[|\varphi_i(\tau)-\varphi_i(\omega)|\leq|(\varphi_i-\varphi)(\tau)-(\varphi_i-\varphi)(\omega)|+|\varphi(\tau)-\varphi(\omega)|\leq ([\varphi_i-\varphi]_\alpha+[\varphi]_\alpha) d_r(\tau,\omega)^\alpha,\]
hence
\[[\varphi_i]_\alpha\leq[\varphi_i-\varphi]_\alpha+[\varphi]_\alpha\]
and as $[\varphi_i-\varphi]_\alpha\to 0$ we get that
\[\sup_i[\varphi_i]_\alpha<\infty\]
which implies that
\[D_*=\sup_i D_{\varphi_i}<\infty.\]

\paragraph{Transfer-operator.}
We recap a few facts about the transfer-operator, for details see e.g. \cite{BSSS2022} and \cite{GiuliettiKloecknerLopesMarcon2018}.
Let $\varphi\in C^\alpha(\Sigma)$ be a Hölder continuous potential on $\Sigma$ ($0<\alpha<1$). Then for $f:\Sigma\to\R$ let
\[\mathcal L_\varphi(f)(x)=\sum_{\sigma y=x}e^{\varphi(y)}f(y)\qquad (x\in\Sigma).\] 
It is known that $\mathcal L_\varphi$ has a positive maximal eigenvalue $\lambda_\varphi$ and a positive (bounded away from 0) eigenfunction $h_\varphi\in C^\alpha(\Sigma)$:
\[\mathcal L_\varphi(h_\varphi)=\lambda_\varphi h_\varphi\]
and the dual operator $\mathcal L^*_\varphi$ has an eigenmeasure $\nu_\varphi$ for eigenvalue $\lambda_\varphi$:
\[\int\mathcal L_\varphi(f)d\nu_\varphi=\lambda_\varphi\int fd\nu_\varphi\qquad \text{for all } f\in C(\Sigma).\]
We call $h_\varphi$ normalized if $\int h_\varphi d\nu_\varphi=1$.
Moreover
\[\log\lambda_\varphi=P(\varphi)\]
where 
\[P(\varphi)=\lim_{n\to\infty}\frac{1}{n}\log\sum_{|\mathbf i|=n}\exp\left(\sup_{j\in[\mathbf i]}\sum_{k=0}^{n-1}\varphi(\sigma^kj)\right),\] 
is the topological pressure function (see e.g. \cite[Remark~7.5]{GiuliettiKloecknerLopesMarcon2018}).
Iterating till $n\in\N$ and rearranging
\begin{equation}\label{eq:transfer-operator-n}
    \int fd\nu_\varphi=e^{-nP(\varphi)}\int\mathcal L^n_\varphi(f)d\nu_\varphi.
\end{equation}
\par The associated shift-invariant Gibbs-measure $\mu_\varphi$ is
\[\frac{d\mu_\varphi}{d\nu_\varphi}=\frac{h_\varphi}{\int h_\varphi d\nu_\varphi}.\]

%--------------------------------------------
\paragraph{Entropy.} We follow the notations of \cite{FengHu09}.
For subsets $\mathcal A, \mathcal C$ of $\mathrm{P}(\Sigma)$, both $\widehat{\mathcal A}$ and $\sigma(\mathcal A)$ will denote the $\sigma$-algebra generated by $\mathcal A$, while $\mathcal A\vee\mathcal C$ will denote $\sigma(\mathcal A\cup\mathcal C)$.

For measurable partitions $\xi, \eta$ of $\Sigma$, $\xi\vee\eta$ will denote their common refinement.
For a finite measurable partition $\xi$ of $\Sigma$ and a $\sigma$-algebra $\mathcal A$ the conditional information of $\xi$ given $\mathcal A$
\[\mathbf{I}_\mu(\xi\mid \mathcal A)=-\sum_{A\in\xi}\mathbbm{1}_A\log \mathbf{E}_\mu(\mathbbm{1}_A\mid \mathcal A).\]
The conditional entropy of $\xi$ given $\mathcal A$
\[H_\mu(\xi\mid \mathcal A)=\int\mathbf{I}_\mu(\xi\mid \mathcal A)d\mu
=-\sum_{A\in \xi}\int_A
\log \mathbf{E}_\mu(\mathbbm{1}_A\mid\mathcal A)\,d\mu.\]
For a finite measurable partition $\eta$, we write
$H_\mu(\xi\mid\eta)=H_\mu(\xi\mid\widehat{\eta})$.

For a topological space $X$, let $\cB(X)$ denote the Borel subsets of $X$. Set $\cB=\cB(\R)$. For measure $\mu$ on $\Sigma$ and IFS $\cF=(f_1,\dots,f_m)$ set
\[h_{\Pi_{\cF}}(\sigma,\mu)=H_\mu(\cP\mid\sigma^{-1}\Pi^{-1}_\cF\cB)-H_\mu(\cP\mid\Pi^{-1}_\cF\cB),\]
where $\cP=\{[i]:1\leq i\leq m\}$. We call $h_{\Pi_{\cF}}(\sigma,\mu)$ the projection entropy.

For $n\geq 1$, let $\cP^n=\{[\mathbf u]:\mathbf u\in\Sigma_n\}$
be the partition into cylinders of length $n$. For a shift-invariant
probability measure $\mu$, the measure-theoretic entropy of the shift is
\[
h_\mu(\sigma)
=\lim_{n\to\infty}-\frac{1}{n}
  \sum_{A\in\cP^n}\mu(A)\log\mu(A)
  =\lim_{n\to\infty}\frac1n
H_\mu\!\left(\bigvee_{k=0}^{n-1}\sigma^{-k}\mathcal P\right),
\]
with the convention $0\log 0=0$.

\paragraph{The Feng--Hu dimension formula.}
Let $\mu$ be a shift-invariant ergodic probability measure on $\Sigma$,
and let $\cF$ be a conformal IFS on the real line. By
\cite[Theorem~2.8]{FengHu09}, the measure
$\nu_{\cF}=(\Pi_{\cF})_*\mu$ is exact-dimensional, with almost
everywhere local dimension $h_{\Pi_{\cF}}(\sigma,\mu)/\chi_\mu(\cF)$.
Consequently,
\begin{equation}\label{eq:feng-hu-dimension}
\dimh\nu_{\cF}
=\frac{h_{\Pi_{\cF}}(\sigma,\mu)}{\chi_\mu(\cF)}.
\end{equation}

\paragraph{The variational principle for pressure.}
The variational principle for pressure (\cite[Theorem~9.10]{Walters82}) states that for a continuous potential ${\varphi:\Sigma\to\R}$
\begin{equation}\label{eq:variational-principle}
P(\varphi)=\sup_{\nu\in\mathcal M_\sigma(\Sigma)}\left(h_\nu(\sigma)+\int_\Sigma\varphi\,d\nu\right),
\end{equation}
where $\mathcal M_\sigma(\Sigma)$ is the set of shift-invariant Borel probability measures on $\Sigma$. For a Hölder continuous potential $\varphi$, the variational principle has a unique maximizer, also called the  equilibrium state, namely the shift-invariant Gibbs measure $\mu_\varphi$ associated with the potential.
Hence by \cite[Theorem~1.22]{Bowen1975},
\begin{equation}\label{eq:vp-maximizer}
P(\varphi)-\int\varphi\,d\mu_{\varphi}=h_{\mu_\varphi}(\sigma).
\end{equation}

%------------------------------------------------------------------------------------------------------------------------------The quasi-Bernoulli case---------------------------
\section{The quasi-Bernoulli case}

%------------------------------------------------------------------------------------------------------------------------------The conformal version of the ``good covering'' lemma---------------------------

As in \cite[Lemma~2.4]{PeresSolomyak00} and before, we define the cut-set for $n\in\N$ when $2^{-n}<|I|$ holds as
\begin{equation}\label{eq:interval-cut-set}
  W_n^{\mathcal F}=W_n=\bigl\{\mathbf i\in\Sigma_*: 
  |\Lambda_{\mathbf{i}}|\le 2^{-n}<|\Lambda_{\mathbf{i}^{-}}|
  \bigr\},
\end{equation}
where $\mathbf i^-$ is obtained from $\mathbf i$ by deleting its last symbol. In order to define $W_n$ for all $n\in\N$, we can have the following conventions:
\[W_n=\{\emptyset\},\qquad f_\emptyset=\mathrm{id},\qquad \Lambda_\emptyset=\Lambda_{\mathcal F}\qquad\text{when }2^{-n}\geq|\Lambda_{\mathcal F}|.\]

We recall the definition from \cite[Definition~2.1]{PeresSolomyak00}.% and also \cite[Definition~3.6.2]{BSS23}.

\begin{df}\label{def:good-cover}
Let $K\subset\R$ be compact. Fix $M>0$ and $n,N\in\N$. A Borel cover $\{G_i\}_{i=1}^k$ of $K$ is called $(M,2^{-n},N)$-good if
\[
  |G_i|\leq M\cdot 2^{-n}\quad\text{for all }i,
\]
and every dyadic cube $Q\in\mathcal D_n$ intersects at most $N$ of the sets $G_i$.
\end{df}

The following lemma was proven in \cite[Lemma~2.4]{PeresSolomyak00}.

\begin{lem}\label{lem:good-cover}
For $\mathcal F\in\Theta_{\gamma_1,\gamma_2}(I)$, there exist $M=M(\mathcal F)>0$ and $N=N(\mathcal F)\in\N$ such that for all $\ell,n\in\N$ and all $\mathbf i\in W_n$, the collection
\[
  \left\{f_{\mathbf i}^{-1}(Q):
  Q\in\mathcal D_{\ell+n},\ Q\cap\Lambda_{\mathbf i}\neq\emptyset\right\}
\]
is an $(M,2^{-\ell},N)$-good cover of $\Lambda$.
\end{lem}

We will need later that $M=M(\mathcal F)$ can be chosen such that it is constant in a $C^{1+\theta}$ neighborhood of $\mathcal F$ when its attractor, $\Lambda_{\mathcal F}$ is not a singleton. In the proof in \cite{PeresSolomyak00}, $M=C^2L_1\lambda_0^{-1}$, where $L_1$  is the bounded distortion constant,
\[\lambda_0=\min\{|f'_j(\xi)|:1\leq j\leq m,\xi\in I\},\]
and $C$ is determined by 
\[C^{-1}||f'_{\mathbf u}||\leq|\Lambda_{\mathbf u}|\leq C||f'_{\mathbf u}||,\]
\[C^{-1}||f'_{\mathbf u}||\cdot|x-y|\leq|f_{\mathbf u}(x)-f_{\mathbf u}(y)|\leq C||f'_{\mathbf u}||\cdot|x-y|\quad (\mathbf u\in\Sigma_*,x,y\in I).
\]
It can be readily seen that all three parameters can be chosen such that they are constant in a suitable neighborhood of $\mathcal F$.

We are going to use the next fact that is from \cite[Lemma~2.3]{PeresSolomyak00}.

\begin{fact}\label{fact:entropy-balls}
\begin{enumerate}
  \item If $\mu$ is a compactly supported probability measure on $\R$, then for all $k,n\in\N$,
  \[
    h_n(\mu)\leq h_{n+k}(\mu)\leq h_n(\mu)+2^{k}.
  \]
  \item There exists a constant $C_1$ such that
  \begin{equation}\label{eq:entropy-ball-comparison}
    \left|
    \int \log\frac{1}{\mu(B(x,2^{-n}))}\,d\mu(x)-h_n(\mu)
    \right|
    \leq C_1
  \end{equation}
  for every compactly supported probability measure $\mu$ on $\R$ and for all $n\in\N$.
\end{enumerate}
\end{fact}

%The following theorem was proven in \cite[Theorem~1.1]{PeresSolomyak00}.
%\begin{thm}\label{thm:entropy-dim-exists}
%For every $\mathcal F\in\mathcal F_m$ and every probability vector $p=(p_1,\dots,p_m)$, the entropy dimension $\dime(\nu_{\mathcal F,p})$ exists.
%\end{thm}

%-------------------------------------------------------------------------------------------------Existence of entropy dimension and lower semicontinuity----------------------------------
\subsection{Lower semicontinuity of the entropy and Hausdorff dimensions of self-conformal measures}

\begin{thm}\label{thm:lsc-hausdorff}
Let $\mu$ be a quasi-Bernoulli measure on $\Sigma$. For $0<\gamma_1<\gamma_2<1$ and $m\geq 2$, the map
\begin{equation}\label{eq:hausdorff-dimension-map}
  \mathcal F\longmapsto \dimh((\Pi_{\mathcal F})_*\mu),
  \qquad \mathcal F\in\Theta_{\gamma_1,\gamma_2}(I)
\end{equation}
is lower semicontinuous.
\end{thm}

\begin{proof}
First we assume that $\mu$ is shift-invariant.

\medskip
\textbf{1. Preliminaries on the quasi-Bernoulli property.}

Let $\mu$ be a quasi-Bernoulli measure on $\Sigma$ with constant $K\geq 1$, i.e. for all $\tbu,\tbv\in\Sigma_*$
\begin{equation}\label{eq:Bernoulli}
K^{-1}\mu([\tbu])\mu([\tbv])\leq\mu([\tbu\tbv])\leq K\mu([\tbu])\mu([\tbv]).
\end{equation}

For $\mu([\tbu])>0$ define the conditional tail measure (see \cite[Definition~1.1]{FraserPollicott2017})
\[\mu_\tbu(A)=\frac{\mu([\tbu]\cap\sigma^{-|\tbu|}(A))}{\mu([\tbu])}.\]

First applying \eqref{eq:Bernoulli} for cylinders, and then extending for Borel sets, we get that
\[K^{-1}\mu\leq\mu_\tbu\leq K\mu.\]

Using
\[\nu_\cF=(\Pi_\cF)_*\mu\qquad\nu_{\cF,\tbu}=(\Pi_\cF)_*\mu_\tbu\]
we have
\[K^{-1}\nu_\cF\leq\nu_{\cF,\tbu}\leq K\nu_\cF.\]

For our cut set $W_n^\cF$
\[\boxed{\nu_\cF=\sum_{\tbu\in W_n^\cF}\mu([\tbu])(f_\tbu)_*\nu_{\cF,\tbu}.}\]
This follows from $\Pi_\cF(\tbu\tau)=f_\tbu(\Pi_\cF(\tau))$.

\medskip
\textbf{2. Estimate on the dyadic entropy.}

Fix $\mathcal F$, and abbreviate $\nu=\nu_{\mathcal F}$. Define
\[\boxed{b_\ell^K(\nu)=\inf\{h_\ell(\eta):\eta\text{ is a probability measure and }K^{-1}\nu\leq\eta\leq K\nu\}.}\]
For $n,\ell\in\N$ we are going to prove
\begin{equation}\label{eq:entropy-superadditivity}
  \boxed{h_{n+\ell}(\nu)\geq h_n(\nu)+b_\ell^K(\nu)-C}
\end{equation}
with a constant $C$.

Let
\[
  \varphi(x)=
  \begin{cases}
    x\log(1/x),&x>0,\\
    0,&x=0.
  \end{cases}
\]
Then
\begin{equation}\label{eq:phi-product}
  \varphi(xy)=y\varphi(x)+x\varphi(y),
  \qquad x,y>0.
\end{equation}
For a Borel set $Q$ and $\tbu\in\Sigma_*$
\[\mu([\tbu]\cap\Pi^{-1}(Q))=\mu([\tbu]\cap\sigma^{-|\tbu|}(\Pi^{-1}(f^{-1}_\tbu(Q))))\]
\[=\mu([\tbu])\mu_\tbu(\Pi^{-1}(f^{-1}_\tbu(Q)))=\mu([\tbu])\nu_{\cF,\tbu}(f^{-1}_\tbu(Q))=\mu([\tbu])(f_\tbu)_*\nu_{\cF,\tbu}(Q).\]
Summing over all $\tbu\in W_n$ we have
\begin{equation}\label{eq:cube-measure-decomposition}
\nu_\cF(Q)=\sum_{\tbu\in W_n(Q)}\mu([\tbu])(f_\tbu)_*\nu_{\cF,\tbu}(Q).
\end{equation}
where
\[
  W_n(Q)=W_n^\cF(Q)=\{\tbu\in W_n: \Lambda_{\tbu}\cap Q\neq\emptyset\}.
\]
Write
\begin{equation}\label{eq:entropy-phi-sum}
  h_k=h_k(\nu)=\sum_{Q\in\mathcal D_k}\varphi(\nu(Q)).
\end{equation}
For $Q\in\mathcal D_{\ell+n}$ define
\[
  p(Q)=\sum_{\mathbf u\in W_n(Q)}\mu([\tbu]).
\]
Clearly
\[p(Q)=\mu\left(\bigcup_{\mathbf u\in W_n(Q)}[\mathbf u]\right)
\leq\mu\left(\Pi^{-1}\left(\bigcup_{\mathbf u\in W_n(Q)}\Lambda_{\mathbf u}\right)\right)
=\nu\left(\bigcup_{\mathbf u\in W_n(Q)}\Lambda_{\mathbf u}\right).\]

If $x\in Q\in\mathcal D_{\ell+n}$, then
\begin{equation}\label{eq:cylinder-weight-ball-bound}
  %\nu(B(x,2^{-n}))\leq
  p(Q)
  \leq \nu(B(x,2^{-n+1})),
\end{equation}
where we used that $|\Lambda_{\mathbf i}|\leq 2^{-n}$ for ${\mathbf i}\in W_n$, and
\[\Lambda_{\mathbf i}\cap Q\ne\emptyset\Longrightarrow \Lambda_{\mathbf i}\subset \bar{B}(x,2^{-n}+2^{-\ell-n})\subset B(x,2^{-n+1}).\]

Apply \eqref{eq:phi-product} with $x=p(Q)$ and $y=\nu(Q)/p(Q)$:
\begin{equation}\label{eq:phi-decomposition}
  \varphi(\nu(Q))=A_1(Q)+A_2(Q),
\end{equation}
where
\begin{equation}\label{eq:first-entropy-term}
  A_1(Q)=\frac{\nu(Q)}{p(Q)}\varphi(p(Q))
  =\nu(Q)\log\frac{1}{p(Q)},
\end{equation}
and, using \eqref{eq:cube-measure-decomposition},
\[
  A_2(Q)=p(Q)\,
  \varphi\left(\frac{\nu(Q)}{p(Q)}\right)
  =p(Q)\,
  \varphi\left(
    \sum_{\mathbf u\in W_n(Q)}
    \frac{\mu([\tbu])}{p(Q)}\nu_{\cF,\tbu}(f_{\mathbf u}^{-1}(Q))
  \right).
\]
Define $A_1(Q)=A_2(Q)=0$ when $p(Q)=0$; also omit the terms where $\mu([\tbu])=0$.

Since $\varphi$ is concave,
\begin{equation}\label{eq:second-entropy-term-bound}
  A_2(Q)
  \geq
  \sum_{\mathbf u\in W_n(Q)}\mu([\tbu])\,
  \varphi\bigl(\nu_{\cF,\tbu}(f_{\mathbf u}^{-1}(Q))\bigr).
\end{equation}
Thus
\begin{equation}\label{eq:entropy-decomposition}
  h_{\ell+n}=
  \sum_{Q\in\mathcal D_{\ell+n}}
  A_1(Q)+A_2(Q).
\end{equation}

From \eqref{eq:cylinder-weight-ball-bound} and Fact~\ref{fact:entropy-balls},
\begin{equation}\label{eq:first-term-ball-estimate}
  \sum_{Q\in\mathcal D_{\ell+n}}A_1(Q)
  \geq
  \int \log\frac{1}{\nu(B(x,2^{-n+1}))}\,d\nu(x)
  \geq h_{n-1}(\nu)-C_1
  \geq h_n(\nu)-C_1-2.
\end{equation}
Hence
\begin{equation}\label{eq:first-term-lower-bound}
  \sum_{Q\in\mathcal D_{\ell+n}}A_1(Q)
  \geq h_n(\nu)-C_2.
\end{equation}
On the other hand, \eqref{eq:second-entropy-term-bound} gives
\begin{equation}\label{eq:second-term-sum-bound}
  \sum_{Q\in\mathcal D_{\ell+n}}A_2(Q)
  \geq
  \sum_{\mathbf u\in W_n}\mu([\tbu])
  \sum_{\substack{Q:Q\in\mathcal D_{\ell+n}\\ Q\cap\Lambda_{\mathbf u}\neq\emptyset}}
  \varphi\bigl(\nu_{\cF,\tbu}(f_{\mathbf u}^{-1}(Q))\bigr).
\end{equation}
By Lemma~\ref{lem:good-cover}, the sets $\{f_{\mathbf u}^{-1}(Q):Q\in\mathcal D_{\ell+n},\ Q\cap\Lambda_{\mathbf u}\neq\emptyset\}$ form a good cover of $\Lambda$. 
Therefore $f_{\mathbf u}^{-1}(Q)\subset B(x,M\cdot 2^{-\ell})$ for any ${\mathbf u}\in W_n$ and $x\in f_{\mathbf u}^{-1}(Q)$.
Using Fact~\ref{fact:entropy-balls} again, for each $\mathbf u\in W_n$,
\[
  \sum_{Q:Q\cap\Lambda_{\mathbf u}\neq\emptyset}
  \varphi\bigl(\nu_{\cF,\tbu}(f_{\mathbf u}^{-1}(Q))\bigr)
  =\sum_{Q:Q\cap \Lambda_{\mathbf u}\ne\emptyset}\nu_{\cF,\tbu}(f_{\mathbf u}^{-1}(Q))\log\frac{1}{\nu_{\cF,\tbu}(f_{\mathbf u}^{-1}(Q))}\]
\[\geq\int\log\frac{1}{\nu_{\cF,\tbu}(B(x,M\cdot 2^{-\ell}))}d\nu_{\cF,\tbu}(x)
  \geq h_\ell(\nu_{\cF,\tbu})-C_1-2^{\widehat k},
\]
where $\widehat k=\widehat k(\mathcal F)=\max\{0,\lceil \log_2 M(\mathcal F)\rceil\}$. Since $\sum_{\mathbf u\in W_n}\mu([\tbu])=1$,
\[
  \sum_{Q\in\mathcal D_{\ell+n}}A_2(Q)
  \geq
  \sum_{\mathbf u\in W_n}\mu([\tbu])
  \sum_{Q:Q\cap\Lambda_{\mathbf u}\neq\emptyset}
  \varphi\bigl(\nu_{\cF,\tbu}(f_{\mathbf u}^{-1}(Q))\bigr)
  \geq
  \sum_{\mathbf u\in W_n}\mu([\tbu])(h_\ell(\nu_{\cF,\tbu})-C_3)
\]
\[  \geq \sum_{\mathbf u\in W_n}\mu([\tbu])(b_\ell^K(\nu)-C_3)\geq b_\ell^K(\nu)-C_3.
\]
Combining this with \eqref{eq:entropy-decomposition} and \eqref{eq:first-term-lower-bound} proves \eqref{eq:entropy-superadditivity}.

According to the note after Lemma \ref{lem:good-cover}, $M(\mathcal F)$ can be chosen to be a constant in a neighborhood of $\mathcal F$, hence the same is true for $C_3$, and also for $C$ in \eqref{eq:entropy-superadditivity}.

\medskip
\textbf{3. $b_\ell^K$ preserves the entropy dimension.}

As $\cF$ is conformal, $\mu$ is shift-invariant and ergodic (because quasi-Bernoulli), the
exact-dimensionality of $\nu_\cF$ follows from the theorem of Feng and Hu \cite[Theorem~2.8]{FengHu09}. Let $d=\dime\nu_\cF=\dimh\nu_\cF$.

Claim.
\[\lim_{\ell\to\infty}\frac{b_\ell^K(\nu)}{\ell\log 2}=d.\]
Set
\[J_\ell(x)=\frac{-\log\nu(\Dl(x))}{\ell\log 2}\]
where $\Dl(x)$ is the half open dyadic interval with length $2^{-\ell}$ containing $x$.

\cite[Theorem B.1]{KaenmakiRajalaSuomala2010} gives that
\[\lim_{\ell\to\infty}J_\ell(x)=d\qquad\text{for almost every }x.\]
We also show that it is a convergence in $L_1$ too, i.e.
\[\int |J_\ell(x)-d| d\nu(x)\longrightarrow 0.\]
As the support of $\nu$ is compact, there are finitely many dyadic intervals with positive measure for a fixed $\ell$. Let $A$ denote the number of such intervals of length 1. Then for every $\ell$, the number of intervals from $\Dl$ with positive measure is at most $A2^{l}$. If $t>1$, then
\[J_\ell(x)>t\Longleftrightarrow\nu(\Dl(x))<2^{-\ell t},\]
hence
\[\nu(J_\ell(x)>t)<A2^{\ell}2^{-\ell t}\leq A2^{-\ell(t-1)}\leq A2^{-(t-1)}.\]
Then using the integral formula for the distribution function, for $T>1$ we get that
\[\int_{\{J_\ell>T\}}J_\ell(x) d\nu(x)=\int_0^\infty \nu\left(\mathbbm{1}_{\{J_\ell>T\}}J_\ell>t\right) dt\]
\[=\int_0^T \nu\left(\mathbbm{1}_{\{J_\ell>T\}}J_\ell>t\right) dt+\int_T^\infty \nu\left(\mathbbm{1}_{\{J_\ell>T\}}J_\ell>t\right) dt\]
\[\leq T\nu(J_\ell>T)+\int_T^\infty A2^{-(t-1)} dt\leq TA2^{-(T-1)}+\frac{1}{\log 2}A2^{-(T-1)}.\]
But the right hand side tends to $0$ as $T\to \infty$, and together with the almost-everywhere convergence, it shows the convergence in $L_1$.

If $K^{-1}\nu\leq\eta\leq K\nu$, then $\eta\ll\nu$, and set 
\[r(x)=\frac{d\eta}{d\nu}.\]
Obviously $K^{-1}\leq r(x)\leq K$ for almost every $x$, and $\int r(x) d\nu(x)=1$. Note that both $\eta$ and $\nu$ are compactly supported, hence finitely many intervals have positive measures. Then
\[h_\ell(\eta)=\sum_{Q\in\Dl}-\eta(Q)\log\eta(Q)=\sum_Q-\eta(Q)\log\nu(Q) - \sum_Q\eta(Q)\log\frac{\eta(Q)}{\nu(Q)}.\]
\[\sum_Q-\eta(Q)\log\nu(Q)=\int -r(x)\log\nu(\Dl(x)) d\nu(x)=\ell\log2\int r(x)J_\ell(x) d\nu(x)\]
\[\sum_Q\eta(Q)\log\frac{\eta(Q)}{\nu(Q)}\leq\sum_Q\eta(Q)\log K\leq \log K.\]
Hence
\[\left\lvert\frac{h_\ell(\eta)}{\ell\log2}-d\right\lvert\leq \left\lvert\int r(x)J_\ell(x) d\nu(x)-d\int r(x) d\nu(x)\right\lvert+\frac{\log K}{\ell\log2}\]
\[\leq\int r(x)|J_\ell(x)-d| d\nu(x)+\frac{\log K}{\ell\log2}\leq\int K|J_\ell(x)-d| d\nu(x)+\frac{\log K}{\ell\log2}.\]
If $\ell\to\infty$, then both terms of the right hand side tend to 0, which shows that 
\begin{equation}\label{eq:blkd}
\lim_{\ell\to\infty}\frac{b_\ell^K(\nu)}{\ell\log 2}=\lim_{\ell\to\infty}\frac{b_\ell^K(\nu)-C}{\ell\log 2}=d.
\end{equation}

\medskip
\textbf{4. We show that $d$ is also the supremum of the previous sequence.}

Fix $n_0\in\N$. If $k\in\N$, then
\[  h_{n_0+k\ell}(\nu)\geq h_{n_0+(k-1)\ell}(\nu)+b_\ell^K(\nu)-C\geq\dots \geq h_{n_0}(\nu)+k(b_\ell^K(\nu)-C).
\]
Dividing both sides by $(n_0+k\ell)\log2$ and letting $k\to\infty$, we get that
\[d=\dime\nu\geq \frac{b_\ell^K(\nu)-C}{\ell\log2}\qquad\text{for every }\ell,\]
hence
\[d\geq\sup_\ell\frac{b_\ell^K(\nu)-C}{\ell\log2}\]
but together with \eqref{eq:blkd} we have
\[d=\sup_\ell\frac{b_\ell^K(\nu)-C}{\ell\log2}.\]

\medskip
\textbf{5. We show that $b_\ell^K(\nu)$ is actually a minimum, i.e.}
\[b_\ell^K(\nu)=\min\{h_\ell(\eta):\eta\text{ is a probability measure and }K^{-1}\nu\leq\eta\leq K\nu\}.\]
If $K^{-1}\nu\leq\eta\leq K\nu$, then 
\[h_\ell(\eta)=\sum_{\substack{Q\in\Dl\\ \nu(Q)>0}}-\eta(Q)\log\eta(Q)\]
which is a finite sum as the support of $\nu$ is bounded. Moreover the values on $Q\in\Dl$ by $\eta$ matters only, and not $\eta$ itself.
Hence if we set $\nu_Q=\nu(Q)$ and
\[H_\ell=\{h_\ell(\eta):\eta\text{ is a probability measure and }K^{-1}\nu\leq\eta\leq K\nu\},\]
then
\[H_\ell=\left\{\sum_{\substack{Q\in\Dl\\ \nu(Q)>0}}-W_Q\log W_Q:\sum_Q W_Q=1,\ \forall Q\ K^{-1}\nu_Q\leq W_Q\leq K\nu_Q\right\}\]
where $W:\Dl\to\R$ is a function taking non-zero values on finitely many intervals only. As the possible $W$ functions constitute a compact set in $\R^m$ for some $m$, then the function $\sum_Q-W_Q\log W_Q$ takes its minimum on that set.

\medskip
\textbf{6. We show that $\nu\mapsto b_\ell^K(\nu)$ is lower semicontinuous.}

Let $\nu_i\Rightarrow \nu$. Let us assume that the supports of all $\nu_i$ and $\nu$ are in a common compact set. Let us also assume that for all $Q\in\Dl$, $\nu(\partial Q)=0$ holds. Then $\nu_i(Q)\to \nu(Q)$.

Choose a subsequence of $(\nu_i)$ such that $\lim b_\ell^K(\nu_i)=\liminf b_\ell^K(\nu_i)$, and for simplicity, denote this new sequence by $(\nu_i)$ too.

For $i\in\N$ let $W_i:\Dl\to\R$ be chosen such that 
\[\sum_Q W_i(Q)=1,\ \forall Q\ K^{-1}\nu_i(Q)\leq W_i(Q)\leq K\nu_i(Q)\text{ and }\sum_Q-W_i(Q)\log W_i(Q)=b_\ell^K(\nu_i).\]

As $(W_i)$ can be considered as a bounded sequence in $\R^m$, then we can chose a convergent subsequence $W_{i_j}\to W$, and $W$ satisfies
\[\sum_Q W(Q)=1,\ \forall Q\ K^{-1}\nu(Q)\leq W(Q)\leq K\nu(Q).\]
But then
\[\sum_Q-W(Q)\log W(Q)\geq b_\ell^K(\nu).\]
As
\[\sum_Q-W_i(Q)\log W_i(Q)\longrightarrow\sum_Q-W(Q)\log W(Q)\]
we get that
\begin{equation}\label{eg:liminfd}
\liminf_ib_\ell^K(\nu_i)\geq b_\ell^K(\nu).
\end{equation}

\medskip
\textbf{7. We prove that $\cF\mapsto\dimh\nu_{\cF}$ is lower semicontinuous.}

If the attractor of $\cF$ is a singleton, then $\dimh\nu_\cF=0$ so the claim is evident. Therefore suppose that $\Lambda_\cF$ is not a singleton.

Let $\cF_i\to\cF$. Then $\nu_{\cF_i}\Rightarrow \nu_\cF$. Then the supports of all $\nu_{\cF_i}$ and $\nu_\cF$ are in a common compact set because IFSs have a common compact domain. It is possible to choose the grid points of $\Dn$ for all $n\in\N$ so that $\forall n\in\N\ \forall\ Q\in\Dn$ $\nu_\cF(\partial Q)=0$ holds. It is because there are countable many atoms of $\nu_\cF$, and countably many grid points for all $\Dn$.

By point (4) for every $\ell\in\N$
\[\liminf_i\dimh\nu_{\cF_i}\geq \liminf_i\frac{b_\ell^K(\nu_{\cF_i})-C}{\ell\log2}\]
holds with a generic constant $C$ (see the discussion at the end of point (2)).
By point (6)
\[\liminf_i\frac{b_\ell^K(\nu_{\cF_i})-C}{\ell\log2}\geq \frac{b_\ell^K(\nu_{\cF})-C}{\ell\log2}.\]
As it holds for every $\ell$ it yields that
\[\liminf_i\dimh\nu_{\cF_i}\geq\sup_\ell\frac{b_\ell^K(\nu_{\cF})-C}{\ell\log2}=\dimh\nu_\cF\]
by point(4).

\medskip
\textbf{8. The theorem in full generality.}

Now we just assume that $\mu$ is quasi-Bernoulli without property shift-invariance. Then \cite[Proposition~5.1 and Corollary~5.3]{Heurteaux2007} give that there is a quasi-Bernoulli and shift-invariant measure $\tilde\mu$ with constant $C\geq 1$ such that
\[C^{-1}\mu\leq\tilde\mu\leq C\mu.\]
For the push-forward measures $\nu_\cF$ and $\tilde\nu_\cF$ we get that
\[C^{-1}\nu_\cF\leq\tilde\nu_\cF\leq C\nu_\cF.\]
Substituting $B(x,r)$ and then taking the logarithm
\[|\log \nu_\cF(B(x,r))-\log \tilde\nu_\cF(B(x,r))|\leq\log C.\]
Dividing both sides with $|\log r|$, and letting $r\to 0$, we get that $\nu_\cF$ is also exact dimensional with the same dimension.
\end{proof}

We can strengthen the previous result in the following way.

\begin{thm}\label{thm:gen}Let $\mu_i$ be quasi-Bernoulli measures on $\Sigma$ with a common quasi-Bernoulli constant such that $\mu_i\Rightarrow\mu$. For $0<\gamma_1<\gamma_2<1$ and $m\geq 2$, 
let $\mathcal F,\mathcal F_i\in\Theta_{\gamma_1,\gamma_2}(I)$ and $\mathcal F_i\to\mathcal F$. Then
\begin{equation}
  \liminf_{i\to\infty}\dimh(\Pi_{\mathcal F_i})_*\mu_i\geq\dimh(\Pi_{\mathcal F})_*\mu.
\end{equation}
\end{thm}
\begin{proof}Set $\nu_i=(\Pi_{\mathcal F_i})_*\mu_i$ and $\nu=(\Pi_{\mathcal F})_*\mu$.

First note that we can assume that all $\mu_i$ is also shift-invariant. In \cite[Corollary~5.3]{Heurteaux2007}, the constant for strongly equivalence depends only the quasi-Bernoulli constant (actually it is the same). Hence if we have a sequence $(\mu_i)$ of quasi-Bernoulli measures with the same quasi-Bernoulli constant, then there are quasi-Bernoulli and shift-invariant measures $\tilde\mu_i$ with the same constant $C\geq 1$ such that 
\[C^{-1}\mu_i\leq\tilde\mu_i\leq C\mu_i.\]
Then choose a weakly convergent subsequence of $(\tilde\mu_i)$, for which $\tilde\mu_{i_j}\to\tilde\mu$. Because cylinders are clopen, $\tilde\mu$ satisfies the above inequality
\[C^{-1}\mu\leq\tilde\mu\leq C\mu,\]
thus it is strongly equivalent with $\mu$, and also shift-invariant, and \cite[Corollary~5.3]{Heurteaux2007} gives that such measure is unique, which yields that every weakly convergent subsequence of $(\tilde\mu_i)$ tends to the same limit, i.e. $\tilde\mu_i\Rightarrow\tilde\mu$ and $\tilde\mu_i, \tilde\mu$ are quasi-Bernoulli with the same constant and shift-invariant as well and 
\[\dimh(\Pi_{\mathcal F_i})_*\tilde\mu_i=\dimh(\Pi_{\mathcal F_i})_*\mu_i\qquad\dimh(\Pi_{\mathcal F})_*\tilde\mu=\dimh(\Pi_{\mathcal F})_*\mu.\]

The uniform convergence of the coding maps and $\mu_i\Rightarrow\mu$ gives that $\nu_i\Rightarrow\nu$. Then point (6) and (7) in the previous proof yields the required inequality because in (7) we just used the weak convergence of the measures.
\end{proof}

Now we show that we can go further in the case of Gibbs measures.

\begin{thm}For $0<\gamma_1<\gamma_2<1$ and $m\geq 2$, let $\mathcal F,\mathcal F_i\in\Theta_{\gamma_1,\gamma_2}(I)$ and $\mathcal F_i\to\mathcal F$. Let $0<\alpha<1$, and $\varphi_i,\varphi\in C^\alpha(\Sigma)$ be Hölder continuous potentials on $\Sigma$ such that $\varphi_i\to\varphi$ in the $C^\alpha$ norm. Let $\mu_{\varphi_i},\mu_\varphi$ denote the (uniquely determined) associated shift-invariant Gibbs measures. Then
\begin{equation}
  \liminf_{i\to\infty}\dimh(\Pi_{\mathcal F_i})_*\mu_{\varphi_i}\geq\dimh(\Pi_{\mathcal F})_*\mu_\varphi.
\end{equation}
\end{thm}
\begin{proof}In order to apply the previous theorem (Theorem \ref{thm:gen}) we need two things: $\mu_{\varphi_i}\Rightarrow\mu_\varphi$, and $\mu_{\varphi_i}$ are quasi-Bernoulli (which is well known) with a common quasi-Bernoulli constant.

\medskip
\textbf{1. $\boldsymbol{\mu_{\varphi_i}}$ are quasi-Bernoulli with a common quasi-Bernoulli constant.}

We are going to prove that there is a common Gibbs constant $G$ which will immediately yield that there is a common quasi-Bernoulli constant as the usual argument gives that $G^3$ can be chosen.

\cite[Corollary B and Corollary 3.7]{GiuliettiKloecknerLopesMarcon2018} gives that the map
\[\mathcal K:C^\alpha(\Sigma)\to C^\alpha(\Sigma),\quad\mathcal K(\varphi)= h_\varphi\]
is analytic in the $C^\alpha(\Sigma)$ norms, where $h_\varphi$ is the normalized leading eigenfunction of the transfer operator $\mathcal L_\varphi$. Then $\mathcal K$ is continuous hence $||\varphi_i-\varphi||_\alpha\to 0$ yields that $||h_{\varphi_i}-h_{\varphi}||_\infty\to 0$. Because $h_\varphi$ is positive, we can choose $a,b\in\R, 0<a<b<\infty$ such that
\[a<h_{\varphi_i}<b\qquad\text{for all }i.\]
Actually this can be seen for all but finitely many functions, and then modifying $a,b$, we can also have it for all functions.

As $h_{\varphi_i}$ is normalized, by
\[\frac{d\mu_{\varphi_i}}{d\nu_{\varphi_i}}=h_{\varphi_i}\]
we get that
\[a\nu_{\varphi_i}([\mathbf u])\leq\mu_{\varphi_i}([\mathbf u])\leq b\nu_{\varphi_i}([\mathbf u])\qquad (\mathbf u\in\Sigma_*).\]

For the transfer operator $\mathcal L_{\varphi_i}$ and for the continuous function $f(x)=\mathbbm{1}_{[\mathbf u]}$, equation \eqref{eq:transfer-operator-n} gives that
\[\nu_{\varphi_i}([\mathbf u])=\int \mathbbm{1}_{[\mathbf u]}d\nu_{\varphi_i}=e^{-nP(\varphi_i)}\int\mathcal L^n_{\varphi_i}(\mathbbm{1}_{[\mathbf u]})d\nu_{\varphi_i}\]
For $n=|\mathbf u|$ we get that
\[\mathcal L^n_{\varphi_i}(\mathbbm{1}_{[\mathbf u]})(x)=\sum_{|\mathbf v|=n}e^{\sum_{k=0}^{n-1}{\varphi_i}(\sigma^k\mathbf vx)}\mathbbm{1}_{[\mathbf u]}(\mathbf vx)\quad (x\in\Sigma),\]
because the induction step is
\[
    \mathcal L^n_{\varphi_i}(\mathbbm{1}_{[\mathbf u]})(x)=\mathcal L_{\varphi_i}(\mathcal L^{n-1}_{\varphi_i}(\mathbbm{1}_{[\mathbf u]}))(x)
    =\sum_{|\mathbf v|=1}e^{{\varphi_i}(\mathbf vx)}\mathcal L^{n-1}_{\varphi_i}(\mathbbm{1}_{[\mathbf u]})(\mathbf vx)
\]
\[
    =\sum_{|\mathbf v|=1}e^{{\varphi_i}(\mathbf vx)}\sum_{|\mathbf w|=n-1}e^{\sum_{k=0}^{n-2}{\varphi_i}(\sigma^k\mathbf w\mathbf vx)}\mathbbm{1}_{[\mathbf u]}(\mathbf w\mathbf vx)
    =\sum_{|\mathbf v|=n}e^{\sum_{k=0}^{n-1}{\varphi_i}(\sigma^k\mathbf vx)}\mathbbm{1}_{[\mathbf u]}(\mathbf vx).
\]
But $\mathbbm{1}_{[\mathbf u]}(\mathbf vx)\ne 0$ iff $\mathbf v=\mathbf u$, hence only one term remains in the sum:
\[\nu_{\varphi_i}([\mathbf u])=e^{-nP(\varphi_i)}\int e^{\sum_{k=0}^{n-1}{\varphi_i}(\sigma^k\mathbf ux)}d\nu_{\varphi_i}(x).\]
Now for $\omega\in[\mathbf u]$
\[\left\lvert\sum_{k=0}^{n-1}{\varphi_i}(\sigma^k\mathbf ux)-\sum_{k=0}^{n-1}{\varphi_i}(\sigma^k\omega)\right\lvert\leq\sum_{k=1}^\infty\var_k(\varphi_i)=D_{\varphi_i}\leq\sup_j D_{\varphi_j}=D_*.\]
Hence
\[e^{-nP(\varphi_i)+\sum_{k=0}^{n-1}{\varphi_i}(\sigma^k\omega)}e^{-D_*}\leq \nu_{\varphi_i}([\mathbf u])\leq e^{-nP(\varphi_i)+\sum_{k=0}^{n-1}{\varphi_i}(\sigma^k\omega)}e^{D_*},\]
which gives that
\[G^{-1}
\leq\frac{\mu_{\varphi_i}([\mathbf u])}{e^{-nP(\varphi_i)+\sum_{k=0}^{n-1}{\varphi_i}(\sigma^k\omega)}}
\leq G,\qquad G=\max\{b,a^{-1}\}e^{D_*}.\]

\medskip
\textbf{2. $\boldsymbol{\mu_{\varphi_i}\Rightarrow\mu_\varphi}$.}

First note that for the topological pressure function $P$,
$\varphi_i\to\varphi$ implies that
\[P(\varphi_i)\to P(\varphi),\]
because
\[|P(\varphi_i)-P(\varphi)|\leq||\varphi_i-\varphi||_\infty.\]
Now choose a convergent subsequence from $(\mu_{\varphi_i})$: $\mu_{\varphi_{i_j}}\Rightarrow\mu$. As $\Sigma$ is compact, it is possible. Set $\mu_j=\mu_{\varphi_{i_j}}$. Now
\[G^{-1}\leq\frac{\mu_j([\mathbf i])}{\exp(-nP(\varphi_{i_j})+\sum_{k=0}^{n-1}\varphi_{i_j}(\sigma^kw))}\leq G\]
holds for all $j\in\N$, all $\mathbf i\in\Sigma_*$, all $w\in[\mathbf i]$, where $|\mathbf i|=n$. By the above observation, for fixed $\mathbf i\in\Sigma_*$ and $w\in[\mathbf i]$, we have
\[\lim_j\frac{\mu_j([\mathbf i])}{\exp(-nP(\varphi_{i_j})+\sum_{k=0}^{n-1}\varphi_{i_j}(\sigma^kw))}=\frac{\mu([\mathbf i])}{\exp(-nP(\varphi)+\sum_{k=0}^{n-1}\varphi(\sigma^kw))}\]
hence
\[G^{-1}\leq\frac{\mu([\mathbf i])}{\exp(-nP(\varphi)+\sum_{k=0}^{n-1}\varphi(\sigma^kw))}\leq G\]
holds. As $\mu_j$ is shift-invariant, so is $\mu.$ Hence $\mu$ is the \textbf{unique} Gibbs measure associated with $\varphi$, i.e. $\mu=\mu_\varphi$ and $\mu_{\varphi_{i_j}}\Rightarrow\mu_\varphi$. As it holds for all convergent subsequences, we get that $\mu_{\varphi_i}\Rightarrow\mu_\varphi$.
\end{proof}

%\begin{rem}Analyzing the previous proof, one can see that the two conditions that $\varphi_i,\varphi\in C^\alpha(\Sigma)$ and $\varphi_i\to\varphi$ in $C^\alpha$ can be weakened the the following.
%\[\varphi_i,\varphi\text{ are Hölder continuous functions},\quad \sup_i D_{\varphi_i}<\infty,\quad ||\varphi_i-\varphi||_\infty\to 0.\]
%For example, it is sufficient that
%\[\varphi_i\in C^{\alpha_i}(\Sigma),\varphi\in C^{\alpha}(\Sigma),\quad \inf_i\alpha_i>0,\quad \sup_i [\varphi_i]_{\alpha_i}<\infty,\quad ||\varphi_i-\varphi||_\infty\to 0.\]
%\end{rem}

Combining the preceding lower semicontinuity result (Theorem \ref{thm:lsc-hausdorff}) with the positivity theorem of Hare and Troscheit yields the following consequence for compact families.

\begin{cor}Let $\mu$ be a quasi-Bernoulli measure with full support on $\Sigma$, and $\mathcal K\subset\Theta_{\gamma_1,\gamma_2}$ be a compact set such that if $\cF\in\mathcal K$, then $\Lambda_\cF$ is not a singleton. Then the map $\cF\mapsto\dimh(\Pi_\cF)_*\mu$ attains a strictly positive minimum on $\mathcal K$.
In particular
\[\inf_{\cF\in\mathcal K}\dimh(\Pi_\cF)_*\mu>0.\]
\end{cor}
\begin{proof}By Theorem \ref{thm:lsc-hausdorff} $\cF\mapsto\dimh(\Pi_\cF)_*\mu$ is lower semicontinuous and therefore attains its minimum
on $\mathcal K$ at some $\cF'\in\mathcal K$.
%\[\exists\cF'\in\mathcal K\text{ such that }\dimh(\Pi_{\cF'})_*\mu=\inf_{\cF\in\mathcal K}\dimh(\Pi_\cF)_*\mu=\min_{\cF\in\mathcal K}\dimh(\Pi_\cF)_*\mu.\]
\cite[Corollary 3.7]{HareTroscheit2021} yields that the lower Assouad dimension of the push-forward of a quasi-Bernoulli measure onto non-degenerate self-conformal sets is positive. Since the lower Assouad dimension of a measure does not exceed its Hausdorff dimension, we obtain $\dimh(\Pi_{\cF'})_*\mu>0$.
\end{proof}

%---------------------------------------------------------------------------------Generic dimension formula for quasi-Bernoulli push-forward measures--------------------------------------------------
\subsection{Generic dimension formula for quasi-Bernoulli push-forward measures on a residual subset of $\Theta_{0,1/2}(I)$}

\begin{thm}\label{thm:residual}
Let $\mu$ be a shift-invariant and quasi-Bernoulli measure on $\Sigma$. For
$\mathcal F\in \Theta_{0,1/2}(I)$, let
\[
        \nu_{\mathcal F}:=\mu\circ \Pi_{\mathcal F}^{-1}
\]
be the push-forward measure on the attractor $\Lambda_{\mathcal F}$, and let
\(h_\mu\) and \(\chi_{\mu}(\mathcal F)\)
denote the corresponding entropy and the Lyapunov exponent.

Then, for a residual subset $H$ of
$\Theta_{0,1/2}(I)$, $\mathcal F\in H$ implies that
\[
        \dimh(\nu_{\mathcal F})
        =
        \min\left\{\frac{h_\mu}{\chi_\mu(\mathcal F)},1\right\}.
\]
\end{thm}

\begin{proof}
Fix $0<\gamma_1<\gamma_2<1/2$. We first prove the result on
$\Theta_{\gamma_1,\gamma_2}(I)$.

Define
\[
        \alpha(\mathcal F):=\dimh(\nu_{\mathcal F})
\]
and
\[
        \beta(\mathcal F)
        :=
        \min\left\{\frac{h_\mu}{\chi_\mu(\mathcal F)},1\right\}.
\]

By Theorem \ref{thm:lsc-hausdorff}, for every $0<\gamma_1<\gamma_2<1/2$, the map
\[
        \mathcal F\mapsto \dimh(\nu_{\mathcal F})
\]
is lower semicontinuous on $\Theta_{\gamma_1,\gamma_2}(I)$, i.e. $\alpha$ is lower semicontinuous on
$\Theta_{\gamma_1,\gamma_2}(I)$. On the other hand, $\beta$ is continuous on
$\Theta_{\gamma_1,\gamma_2}(I)$, since
$\mathcal F\mapsto \chi_\mu(\mathcal F)$ is continuous by \cite[Lemma~14.4.3(ii)]{BSS23}. Moreover, the general
upper bound for conformal IFS measures (see \cite[Theorem~14.2.3]{BSS23}) gives
\[
        \alpha(\mathcal F)\leq \beta(\mathcal F)
        \qquad
        \text{for every }
        \mathcal F\in\Theta_{\gamma_1,\gamma_2}(I).
\]

We claim that
\[
\begin{split}
        \{\mathcal F\in\Theta_{\gamma_1,\gamma_2}(I):
        \alpha(\mathcal F)=\beta(\mathcal F)\}
        =
        \{\mathcal F\in\Theta_{\gamma_1,\gamma_2}(I):
        \alpha \text{ is continuous at } \mathcal F\}.
\end{split}
\]

First suppose that $\alpha(\mathcal F)=\beta(\mathcal F)$. If
$\mathcal F_n\to\mathcal F$, then lower semicontinuity gives
\[
        \liminf_{n\to\infty}\alpha(\mathcal F_n)
        \geq
        \alpha(\mathcal F)
        =
        \beta(\mathcal F).
\]
Since $\alpha\leq\beta$ and $\beta$ is continuous, we also have
\[
        \limsup_{n\to\infty}\alpha(\mathcal F_n)
        \leq
        \limsup_{n\to\infty}\beta(\mathcal F_n)
        =
        \beta(\mathcal F)
        =
        \alpha(\mathcal F).
\]
Hence $\alpha(\mathcal F_n)\to\alpha(\mathcal F)$, so $\alpha$ is continuous at
$\mathcal F$.

Conversely, suppose that $\alpha$ is continuous at $\mathcal F=(f_1,\ldots,f_m)$, but
\[
        \alpha(\mathcal F)<\beta(\mathcal F).
\]
Since $\alpha$ is continuous at $\mathcal F$ and $\beta$ is continuous, there is
a neighbourhood $\mathcal U$ of $\mathcal F$ such that
\[
        \alpha(\mathcal G)<\beta(\mathcal G)
        \qquad
        \text{for every } \mathcal G\in\mathcal U.
\]
Consider the translation family
\[
        \mathcal F^t
        :=
        (f_1+t_1,\ldots,f_m+t_m).
\]
Since $\gamma_2<1/2$, the transversality argument (\cite[Section~14.4]{BSS23}) used for translation
families applies, and therefore the measure-theoretic transversality theorems (\cite[Theorem~14.4.2(i)]{BSS23} and \cite[Theorem~14.5.2]{BSS23})
imply that for Lebesgue-a.e. sufficiently small $t$,
\[
        \dimh(\nu_{\mathcal F^t})
        =
        \min\left\{
        \frac{h_\mu}{\chi_\mu(\mathcal F^t)},1
        \right\}.
\]
Equivalently,
\[
        \alpha(\mathcal F^t)=\beta(\mathcal F^t)
\]
for Lebesgue-a.e. sufficiently small $t$. Taking such a small $t$ with
$\mathcal F^t\in\mathcal U$ gives a contradiction to
$\alpha(\mathcal G)<\beta(\mathcal G)$ on $\mathcal U$. Therefore
$\alpha(\mathcal F)=\beta(\mathcal F)$.

Thus the equality set is precisely the set of continuity points of $\alpha$.
The set of continuity points of a real-valued function on a metric space is a
$G_\delta$ set. Moreover, the preceding translation-family argument shows that
the equality set is dense in $\Theta_{\gamma_1,\gamma_2}(I)$. Hence
\[
        \left\{\mathcal F\in\Theta_{\gamma_1,\gamma_2}(I):
        \dimh(\nu_{\mathcal F})
        =
        \min\left\{
        \frac{h_\mu}{\chi_\mu(\mathcal F)},1
        \right\}
        \right\}
\]
is residual in $\Theta_{\gamma_1,\gamma_2}(I)$.

Finally, use the exhaustion
\[
        \Theta_{0,1/2}(I)
        =
        \bigcup_{n\geq 5}
        \Theta_{1/n,\,1/2-1/n}(I).
\]
Since the exhaustion pieces are relatively closed, and the above residuality statement holds on each piece 
$\Theta_{1/n,\,1/2-1/n}(I)$, it follows that the desired equality holds for a
residual subset of $\Theta_{0,1/2}(I)$.
\end{proof}

%-----------------------------------------------------------------------------------------------------------------------------------
\section{The ergodic case}

We prove the lower semicontinuity of the maps $\mathcal F\mapsto h_{\Pi_{\cF}}(\sigma,\mu)$ and ${\mathcal F\mapsto \dimh((\Pi_{\mathcal F})_*\mu)}$ when $\mathcal F\in\Theta_{\gamma_1,\gamma_2}(I)$.

\begin{thm}\label{thm:lsc-ergodic}
Let $\mu$ be a shift-invariant ergodic measure on $\Sigma$. For $0<\gamma_1<\gamma_2<1$ and $m\geq 2$. Let $\mathcal F,\mathcal F_i\in\Theta_{\gamma_1,\gamma_2}(I)$ and $\mathcal F_i\to\mathcal F$ in the $C^1$ norm. Then
\begin{equation}\label{eq:hausdorff-dimension-map1}
  \liminf_{i\to\infty}h_{\Pi_{\cF_i}}(\sigma,\mu)\geq h_{\Pi_{\cF}}(\sigma,\mu).
\end{equation}
In particular, 
\begin{equation}\label{eq:hausdorff-dimension-map2}
  \liminf_{i\to\infty}\dimh(\Pi_{\mathcal F_i})_*\mu\geq\dimh(\Pi_{\mathcal F})_*\mu.
\end{equation}
\end{thm}

\begin{proof}Let $\cB$ denote the Borel subsets of $\R$. Let
\[\cP^n=\{[\mathbf u]:\mathbf u\in\Sigma_n\}\]
be the partition of $\Sigma$ into length-$n$ cylinders. For $\mathcal F\in\Theta_{\gamma_1,\gamma_2}(I), n\in\N$ set
\[a_n(\cF)=H_\mu(\cP^n\mid \Pi_\cF^{-1}\cB).\]

\medskip
\textbf{1. The sequence $(a_n(\cF))$ is subadditive.}

We are going to verify the following calculation:
\[
\begin{aligned}
a_{n+m}(\cF)
&=H_\mu(\cP^{n+m}\mid \Pi_\cF^{-1}\cB)\\
&=H_\mu(\cP^n\vee\sigma^{-n}\cP^m\mid \Pi_\cF^{-1}\cB)\\
&=H_\mu(\cP^n\mid \Pi_\cF^{-1}\cB) + H_\mu(\sigma^{-n}\cP^m\mid \widehat{\cP^n}\vee\Pi_\cF^{-1}\cB)\\
&=H_\mu(\cP^n\mid \Pi_\cF^{-1}\cB) + H_\mu(\sigma^{-n}\cP^m\mid \widehat{\cP^n}\vee\sigma^{-n}\Pi_\cF^{-1}\cB)\\
&\leq a_n(\cF) + H_\mu(\sigma^{-n}\cP^m\mid \sigma^{-n}\Pi_\cF^{-1}\cB)\\
&=a_n(\cF) + H_\mu(\cP^m\mid \Pi_\cF^{-1}\cB)\\
&=a_n(\cF)+a_m(\cF).
\end{aligned}
\]
Clearly $\cP^{n+m}=\cP^n\vee\sigma^{-n}\cP^m$ holds. 

The next step 
\[H_\mu(\cP^n\vee\sigma^{-n}\cP^m\mid \Pi_\cF^{-1}\cB)\\
=H_\mu(\cP^n\mid \Pi_\cF^{-1}\cB) + H_\mu(\sigma^{-n}\cP^m\mid \widehat{\cP^n}\vee\Pi_\cF^{-1}\cB)\]
follows from \cite[Lemma~4.4(iii)]{FengHu09} with
\[\xi=\cP^n,\quad \eta=\sigma^{-n}\cP^m,\quad\mathcal A=\Pi_\cF^{-1}\cB.\]

\cite[Lemma~4.7]{FengHu09} gives that
\[\widehat{\cP^n}\vee\Pi_\cF^{-1}\cB=\widehat{\cP^n}\vee\sigma^{-1}\Pi_\cF^{-1}\cB.\]
Iterating we get that
\[\widehat{\cP^n}\vee\Pi_\cF^{-1}\cB=\widehat{\cP^n}\vee\sigma^{-n}\Pi_\cF^{-1}\cB.\]

Removing some conditioning increases entropy, i.e. 
\[\sigma^{-n}\Pi_\cF^{-1}\cB\subseteq \widehat{\cP^n}\vee\sigma^{-n}\Pi_\cF^{-1}\cB\]
gives
\[H_\mu(\sigma^{-n}\cP^m\mid \widehat{\cP^n}\vee\sigma^{-n}\Pi_\cF^{-1}\cB)
\leq H_\mu(\sigma^{-n}\cP^m\mid \sigma^{-n}\Pi_\cF^{-1}\cB).\]

And lastly
\[H_\mu(\sigma^{-n}\cP^m\mid \sigma^{-n}\Pi_\cF^{-1}\cB) = H_\mu(\cP^m\mid \Pi_\cF^{-1}\cB)\]
holds because $\mu$ is $\sigma$-invariant. (This also follows from \cite[Lemma~4.4(i)]{FengHu09}.)

\medskip
\textbf{2. Consequence of subadditivity}

Fekete's lemma \cite[Theorem~4.10]{Walters82} gives that
\begin{equation}\label{eq:subadd}
  \lim_{n\to\infty}\frac{a_n(\cF)}{n}=\inf_{n}\frac{a_n(\cF)}{n}.
\end{equation}

\medskip
\textbf{3. Supremum formula for projection entropy}

For ergodic $\mu$, \cite[Proposition~4.10]{FengHu09} gives that
\[\lim_{n\to\infty}\frac{1}{n}\mathbf{I}_\mu(\cP^n\mid \Pi_\cF^{-1}\cB)(x)=h_\mu(\sigma)-h_{\Pi_{\cF}}(\sigma,\mu)\]
$\mu$-almost everywhere and in $L^1$. Integrating both sides we get that
\[\lim_{n\to\infty}\frac{1}{n}H_\mu(\cP^n\mid \Pi_\cF^{-1}\cB)=\lim_{n\to\infty}\frac{a_n(\cF)}{n}=h_\mu(\sigma)-h_{\Pi_{\cF}}(\sigma,\mu).\]
Rearranging and using \eqref{eq:subadd}
\begin{equation}\label{eq:projent}
h_{\Pi_{\cF}}(\sigma,\mu)=h_\mu(\sigma)-\inf_{n}\frac{a_n(\cF)}{n}=\sup_{n}\left(h_\mu(\sigma)-\frac{a_n(\cF)}{n}\right).
\end{equation}

\medskip
\textbf{4. $\boldsymbol{a_n(\cF)}$ is upper semicontinuous}

Let $\cF_i\to\cF$ in $C^1$. Then $\Pi_{\cF_i}\to\Pi_{\cF}$ uniformly. Set $\nu=(\Pi_{\cF})_*\mu$.

For $k\in\N$ let $\cQ_k$ be a translation of the $k$-dyadic grid such that ${\nu(\partial_IQ)=0}$ when $Q\in\cQ_k$. As the set of all such boundary points for all $k$ is countable, we can choose a common such translation. We just keep those intervals in $\cQ_k$ which intersect $I$. As $I$ is compact, $\cQ_k$ is finite.

We show that for a fixed $k$
\[\lim_{i\to\infty}H_\mu(\cP^n\mid \Pi_{\cF_i}^{-1}\cQ_k)= H_\mu(\cP^n\mid \Pi_{\cF}^{-1}\cQ_k).\]

%By the definition of conditional entropy
By the entropy rule \cite[Lemma~4.4(iii)]{FengHu09},
\[H_\mu(\xi\mid\eta)=H_\mu(\xi\vee\eta)-H_\mu(\eta)\]
for finite measurable partitions $\xi,\eta$. Hence
\[H_\mu(\cP^n\mid \Pi_{\cF_i}^{-1}\cQ_k)=H_\mu(\cP^n\vee \Pi_{\cF_i}^{-1}\cQ_k) - H_\mu(\Pi_{\cF_i}^{-1}\cQ_k)\]
%\sum_{j=1}^{r_k}\mu(\Pi_{\cF_i}^{-1}Q_j)H_\mu(\cP^n\mid \Pi_{\cF_i}^{-1}Q_j))\]
%\[=\sum_{j=1}^{r_k}\mu(\Pi_{\cF_i}^{-1}Q_j)\sum_{\mathbf u\in\Sigma_n} -\mu([\mathbf u]\mid \Pi_{\cF_i}^{-1}Q_j)\log\mu([\mathbf u]\mid \Pi_{\cF_i}^{-1}Q_j)\]
%\[=\sum_{j=1}^{r_k}\mu(\Pi_{\cF_i}^{-1}Q_j)\sum_{\mathbf u\in\Sigma_n} -\frac{\mu([\mathbf u]\cap \Pi_{\cF_i}^{-1}Q_j)}{\mu(\Pi_{\cF_i}^{-1}Q_j)}\log\frac{\mu([\mathbf u]\cap \Pi_{\cF_i}^{-1}Q_j)}{\mu(\Pi_{\cF_i}^{-1}Q_j)}\]
%\[=\sum_{j=1}^{r_k}\sum_{\mathbf u\in\Sigma_n} -\mu([\mathbf u]\cap \Pi_{\cF_i}^{-1}Q_j)\log\mu([\mathbf u]\cap \Pi_{\cF_i}^{-1}Q_j) + \sum_{j=1}^{r_k}\log\mu(\Pi_{\cF_i}^{-1}Q_j)\sum_{\mathbf u\in\Sigma_n}\mu([\mathbf u]\cap \Pi_{\cF_i}^{-1}Q_j)\]
%\[=\sum_{j=1}^{r_k}\sum_{\mathbf u\in\Sigma_n} -\mu([\mathbf u]\cap \Pi_{\cF_i}^{-1}Q_j)\log\mu([\mathbf u]\cap \Pi_{\cF_i}^{-1}Q_j) + \sum_{j=1}^{r_k}\mu(\Pi_{\cF_i}^{-1}Q_j)\log\mu(\Pi_{\cF_i}^{-1}Q_j)\]
\begin{equation}\label{eq:ssp}
=\sum_{Q\in\cQ_k}\sum_{\mathbf u\in\Sigma_n} \varphi(\mu([\mathbf u]\cap \Pi_{\cF_i}^{-1}Q)) - \sum_{Q\in\cQ_k}\varphi(\mu(\Pi_{\cF_i}^{-1}Q))
\end{equation}
where $\varphi(x)=-x\log x$ and $\varphi(0)=0$.

If $\mathbf u\in\Sigma_n$ and $Q\in\cQ_k$
\begin{equation}\label{eq:muint}
\mu([\mathbf u]\cap \Pi_{\cF_i}^{-1}Q)=\int_\Sigma\mathbbm{1}_{[\mathbf u]}(\omega)\mathbbm{1}_{Q}(\Pi_{\cF_i}(\omega)) d\mu(\omega).
\end{equation}
If for $\omega\in\Sigma$, 
\[\Pi_\cF(\omega)\notin\partial_IQ,\]
then $\omega$ belongs to a set of full measure because
\[\mu(\Pi_\cF^{-1}(\partial_IQ))=\nu(\partial_IQ)=0.\]
If $\Pi_\cF(\omega)$ is not on the boundary, then there is a neighborhood which is a subset of either $\mathrm{int} Q$ or $\mathrm{ext} Q$. As
\[\Pi_{\cF_i}(\omega)\to\Pi_{\cF}(\omega)\]
the points $\Pi_{\cF_i}(\omega)$ eventually will belong to that neighborhood. Hence
\[\mathbbm{1}_{Q}(\Pi_{\cF_i}(\omega))\to\mathbbm{1}_{Q}(\Pi_{\cF}(\omega))\quad\text{for }\mu\text{-almost every }\omega.\]
As both functions are bounded by $1$ in \eqref{eq:muint}, we get that
\[\lim_i\mu([\mathbf u]\cap \Pi_{\cF_i}^{-1}Q)=\mu([\mathbf u]\cap \Pi_{\cF}^{-1}Q).\]
Similarly we can get that
\[\lim_i\mu(\Pi_{\cF_i}^{-1}Q)=\mu(\Pi_{\cF}^{-1}Q).\]
Then applying these in \eqref{eq:ssp}, we can take the limit because $\varphi$ is continuous and there are finitely many terms in the sum. Therefore we get that
\[\lim_i\left(\sum_{Q\in\cQ_k}\sum_{\mathbf u\in\Sigma_n} \varphi(\mu([\mathbf u]\cap \Pi_{\cF_i}^{-1}Q)) - \sum_{Q\in\cQ_k}\varphi(\mu(\Pi_{\cF_i}^{-1}Q))\right)\]
\[=\sum_{Q\in\cQ_k}\sum_{\mathbf u\in\Sigma_n} \varphi(\mu([\mathbf u]\cap \Pi_{\cF}^{-1}Q)) - \sum_{Q\in\cQ_k}\varphi(\mu(\Pi_{\cF}^{-1}Q)).\]
In other words
\[\lim_iH_\mu(\cP^n\mid \Pi_{\cF_i}^{-1}\cQ_k)=H_\mu(\cP^n\mid \Pi_{\cF}^{-1}\cQ_k).\]
As $\Pi_{\cF_i}^{-1}\cQ_k\subseteq\Pi_{\cF_i}^{-1}\cB$,
\begin{equation}\label{eq:afH}
a_n(\cF_i)\leq H_\mu(\cP^n\mid \Pi_{\cF_i}^{-1}\cQ_k),
\end{equation}
holds which gives that for all $k$
\[\limsup_ia_n(\cF_i)\leq H_\mu(\cP^n\mid \Pi_{\cF}^{-1}\cQ_k).\]
Using
\[\sigma(\Pi_{\cF}^{-1}\cQ_k)\subseteq \sigma(\Pi_{\cF}^{-1}\cQ_{k+1}),\qquad\sigma\left(\bigcup_k\Pi_{\cF}^{-1}\cQ_k\right)=\Pi_{\cF}^{-1}\cB,\]
\cite[Lemma~4.4(iv)]{FengHu09} gives that
\[\lim_k H_\mu(\cP^n\mid \Pi_{\cF}^{-1}\cQ_k)=H_\mu(\cP^n\mid \Pi_\cF^{-1}\cB).\]
Hence
\[\limsup_ia_n(\cF_i)\leq H_\mu(\cP^n\mid \Pi_\cF^{-1}\cB)=a_n(\cF).\]
holds.

\medskip
\textbf{5. $\boldsymbol{\mathcal F\mapsto h_{\Pi_{\cF}}(\sigma,\mu)}$ is lower semicontinuous.}

As the supremum of lower semicontinuous functions is lower semicontinuous, \eqref{eq:projent} and the upper semicontinuity of $(a_n(\cF))$ give the claim.

\medskip
\textbf{6. Equation \eqref{eq:hausdorff-dimension-map1} implies equation \eqref{eq:hausdorff-dimension-map2}.}
 
Actually we show that they are equivalent.

If $\mu$ is ergodic, then Feng-Hu's \cite[Theorem~2.8]{FengHu09} gives that $(\Pi_{\mathcal F})_*\mu$ is exact-dimensional with almost everywhere local dimension $h_{\Pi_{\cF}}(\sigma,\mu)/\chi_\mu(\cF)$. By \cite{Young82}, the Hausdorff-dimension of an exact-dimensional measure equals to its almost everywhere local dimension, hence we get
\begin{equation}\label{eq:FHm}
\dimh((\Pi_{\mathcal F})_*\mu)=\frac{h_{\Pi_{\cF}}(\sigma,\mu)}{\chi_\mu(\cF)},
\end{equation}
where 
\[\chi_\mu(\mathcal F)
 =-\int_\Sigma
 \log\left|f_{\omega_1}'\bigl(\Pi_{\mathcal F}(\sigma\omega)\bigr)\right|
 \,d\mu(\omega).\]
We show that
\begin{equation}\label{eq:Lyc}
\mathcal F\longmapsto \chi_\mu(\cF)\quad\text{is continuous}
\end{equation}
in the $C^1$ topology.

Let $\cF_j=(f_{j,1},\dots,f_{j,m}), \cF=(f_1,\dots,f_m)$ and $\cF_j\to\cF$ in the $C^1$ topology i.e.
\[||f_{j,i}-f_i||_\infty+||f'_{j,i}-f'_i||_\infty\to 0\quad\text{for all }i.\]
This yields that 
\[||\Pi_{\cF_j}-\Pi_{\cF}||_\infty\to 0.\]
Therefore
\[f'_{j,\omega_1}\bigl(\Pi_{\cF_j}(\sigma\omega)\bigr)\longrightarrow f'_{\omega_1}\bigl(\Pi_{\cF}(\sigma\omega)\bigr)\quad\text{uniformly on }\Sigma.\]
As $\mathcal F_j,\mathcal F\in\Theta_{\gamma_1,\gamma_2}(I)$, we have $0<\gamma_1\leq|f'_{j,\omega_1}|, |f'_{\omega_1}|\leq\gamma_2$ which gives that
\[\log\left|f'_{j,\omega_1}\bigl(\Pi_{\cF_j}(\sigma\omega)\bigr)\right|\longrightarrow \log\left|f'_{\omega_1}\bigl(\Pi_{\cF}(\sigma\omega)\bigr)\right|\quad\text{uniformly on }\Sigma,\]
which implies \eqref{eq:Lyc}.

Now \eqref{eq:FHm} and \eqref{eq:Lyc} give that
\[\mathcal F\longmapsto h_{\Pi_{\cF}}(\sigma,\mu)\]
is lower semicontinuous if and only if
\[\mathcal F\longmapsto \dimh((\Pi_{\mathcal F})_*\mu)\]
is lower semicontinuous.
\end{proof}

\begin{rem}Despite the functions are from $\Theta_{\gamma_1,\gamma_2}(I)$, in the theorem we just needed $C^1$ convergence, not $C^{1+\theta}$.
\end{rem}

Now we are going to answer the following question: how much can the projected dimension drop when the symbolic measure also varies?

\begin{thm}\label{thm:lsc-joint}
Let $\mu,\mu_i$ be shift-invariant ergodic measures on $\Sigma$ such that $\mu_i\Rightarrow\mu$. For $0<\gamma_1<\gamma_2<1$ and $m\geq 2$. Let $\mathcal F,\mathcal F_i\in\Theta_{\gamma_1,\gamma_2}(I)$ and $\mathcal F_i\to\mathcal F$ in the $C^1$ norm. Set
\[d_i=\dimh(\Pi_{\mathcal F_i})_*\mu_i,\qquad d=\dimh(\Pi_{\mathcal F})_*\mu.\]
Then
\begin{equation}\label{eq:hausdorff-dimension-map3}
  \liminf_{i\to\infty}d_i\geq d-\frac{h_\mu(\sigma)-\liminf_{i\to\infty}h_{\mu_i}(\sigma)}{\chi_\mu(\cF)}.
\end{equation}
\end{thm}

\begin{proof}We indicate the changes to the previous proof needed when the symbolic measure varies.

Define
\[a_n(\cF,\eta)=H_\eta(\cP^n\mid \Pi_\cF^{-1}\cB),\]
which is subadditive by point 1. Hence we have the formula
\[h_\eta(\sigma)-h_{\Pi_{\cF}}(\sigma,\eta)=\inf_{n}\frac{a_n(\cF,\eta)}{n}.
\]
Modifying point 2, we can show that for a fixed $n$, $a_n(\cF,\eta)$ is jointly upper semicontinuous. In order to show 
\[\lim_{i\to\infty}H_{\mu_i}(\cP^n\mid \Pi_{\cF_i}^{-1}\cQ_k)= H_{\mu}(\cP^n\mid \Pi_{\cF}^{-1}\cQ_k),\]
we would need
\[\lim_i\mu_i([\mathbf u]\cap \Pi_{\cF_i}^{-1}Q)=\mu([\mathbf u]\cap \Pi_{\cF}^{-1}Q) \quad\text{and}\quad
\lim_i\mu_i(\Pi_{\cF_i}^{-1}Q)=\mu(\Pi_{\cF}^{-1}Q).\]
Set \(\Pi_i=\Pi_{\mathcal F_i}\), \(\Pi=\Pi_{\mathcal F}\), and \(C=[\mathbf u]\). Clearly
\[|\mu_i(C\cap \Pi_i^{-1}Q)-\mu(C\cap \Pi^{-1}Q)|\leq\]
\begin{equation}\label{eq:uQ}
|\mu_i(C\cap \Pi_i^{-1}Q)-\mu_i(C\cap \Pi^{-1}Q)| + |\mu_i(C\cap \Pi^{-1}Q)-\mu(C\cap \Pi^{-1}Q)|.
\end{equation}
Regarding the first term,
\[\left|\mu_i\bigl(C\cap\Pi_i^{-1}(Q)\bigr)
      -\mu_i\bigl(C\cap\Pi^{-1}(Q)\bigr)\right|
\]
\[= \left|\mu_i\bigl(C\cap(\Pi_i^{-1}(Q)-\Pi^{-1}(Q))\bigr)
      -\mu_i\bigl(C\cap(\Pi^{-1}(Q)-\Pi_i^{-1}(Q))\bigr)\right|
\]
\[
\le \mu_i\bigl(\Pi_i^{-1}(Q)-\Pi^{-1}(Q)\bigr)+\mu_i\bigl(\Pi^{-1}(Q)-\Pi_i^{-1}(Q)\bigr).
\]
Now fix $\delta>0$. As $\Pi_i\to\Pi$ uniformly, for all sufficiently large $i$, ${|\Pi_i(\omega)-\Pi(\omega)|<\delta}$. If exactly one of $\Pi_i(\omega)$ and $\Pi(\omega)$ belongs to $Q$, then ${\mathrm{dist}(\Pi(\omega),\partial Q)\leq\delta}$.
Set $F_\delta=\{\omega:\mathrm{dist}(\Pi(\omega),\partial Q)\leq\delta\}.$ Therefore we get
\[\left|\mu_i\bigl(C\cap\Pi_i^{-1}(Q)\bigr)
      -\mu_i\bigl(C\cap\Pi^{-1}(Q)\bigr)\right|\leq\mu_i(F_\delta).\]
As $F_\delta$ is closed, by the closed-set inequality for weak convergence,
\[\limsup_i\mu_i(F_\delta)\leq\mu(F_\delta).\]
When \(\delta\downarrow0\), the set $F_\delta$ decreases to $\Pi^{-1}(\partial Q)$, thus $\mu(F_\delta)\to\mu(\Pi^{-1}(\partial Q))=0$ which gives that the first term in \eqref{eq:uQ} converges to $0$.

Regarding the second term in \eqref{eq:uQ}, as $C$ is clopen,
\[\partial(C\cap \Pi_{\cF}^{-1}Q)\subseteq C\cap \partial(\Pi_{\cF}^{-1}Q)\subseteq C\cap \Pi_{\cF}^{-1}(\partial Q)\subseteq \Pi_{\cF}^{-1}(\partial Q).\]
Then $\mu(\Pi_{\cF}^{-1}(\partial Q))=\nu(\partial Q)=0$ and $\mu_i\Rightarrow\mu$ give that the second term in \eqref{eq:uQ} also converges to $0$.

Similarly to \eqref{eq:afH},
\[a_n(\cF_i,\mu_i)\leq H_{\mu_i}(\cP^n\mid \Pi_{\cF_i}^{-1}\cQ_k),\]
holds which gives that for all $k$
\[\limsup_ia_n(\cF_i,\mu_i)\leq H_\mu(\cP^n\mid \Pi_{\cF}^{-1}\cQ_k).\]
Again
\[\lim_k H_\mu(\cP^n\mid \Pi_{\cF}^{-1}\cQ_k)=H_\mu(\cP^n\mid \Pi_\cF^{-1}\cB)\]
yields that
\[\limsup_ia_n(\cF_i,\mu_i)\leq H_\mu(\cP^n\mid \Pi_\cF^{-1}\cB)=a_n(\cF,\mu).\]

Therefore we get that
\[h_\eta(\sigma)-h_{\Pi_{\cF}}(\sigma,\eta)\]
is the infimum of upper semicontinuous functions, hence it is upper semicontinuous.

Consider
\[\liminf_ih_{\Pi_{\cF_i}}(\sigma,\mu_i)=\liminf_i \bigl(h_{\mu_i}(\sigma)-(h_{\mu_i}(\sigma)-h_{\Pi_{\cF_i}}(\sigma,\mu_i))\bigr)\]
\[\geq\liminf_i h_{\mu_i}(\sigma)-\limsup_i(h_{\mu_i}(\sigma)-h_{\Pi_{\cF_i}}(\sigma,\mu_i))\]
\[\geq\liminf_i h_{\mu_i}(\sigma)-(h_{\mu}(\sigma)-h_{\Pi_{\cF}}(\sigma,\mu))\]
\[=h_{\Pi_{\cF}}(\sigma,\mu)-(h_{\mu}(\sigma)-\liminf_i h_{\mu_i}(\sigma)).\]
Dividing by $\chi_\mu(\cF)$
\[\liminf_i\frac{h_{\Pi_{\cF_i}}(\sigma,\mu_i)}{\chi_\mu(\cF)}
\geq \frac{h_{\Pi_{\cF}}(\sigma,\mu)}{\chi_\mu(\cF)} 
- \frac{h_\mu(\sigma)-\liminf_i h_{\mu_i}(\sigma)}{\chi_\mu(\cF)}.
\]
For showing 
\[\liminf_i\frac{h_{\Pi_{\cF_i}}(\sigma,\mu_i)}{\chi_\mu(\cF)}=\liminf_i\frac{h_{\Pi_{\cF_i}}(\sigma,\mu_i)}{\chi_{\mu_i}(\cF_i)}=\liminf_i\dimh(\Pi_{\cF_i})_*\mu_i\]
we would need that $\chi_{\mu_i}(\cF_i)\to\chi_\mu(\cF)>0$. To see this,
\[|\chi_{\mu_i}(\cF_i)-\chi_\mu(\cF)|\leq|\chi_{\mu_i}(\cF_i)-\chi_{\mu_i}(\cF)|+|\chi_{\mu_i}(\cF)-\chi_\mu(\cF)|.\]
The first term converges to $0$ because the logarithmic derivative functions converge uniformly by $C^1$ convergence, and so does the second term because $\mu_i\Rightarrow\mu$.
\end{proof}

Because of the upper semicontinuity of the symbolic entropy (\cite[Theorem 8.2]{Walters82}), the two conditions are equivalent: 
\[h_\mu(\sigma)-\liminf_{i\to\infty}h_{\mu_i}(\sigma)\leq 0\quad \Longleftrightarrow \quad \lim_{i\to\infty}h_{\mu_i}(\sigma)= h_\mu(\sigma).\]
Therefore we obtain the following.

\begin{cor}
Let $\mu,\mu_i$ be shift-invariant ergodic measures on $\Sigma$ such that
\[\mu_i\Rightarrow\mu\qquad \text{ and }\qquad  h_{\mu_i}(\sigma)\to h_\mu(\sigma).\]
For $0<\gamma_1<\gamma_2<1$ and $m\geq 2$. Let $\mathcal F,\mathcal F_i\in\Theta_{\gamma_1,\gamma_2}(I)$ and $\mathcal F_i\to\mathcal F$ in the $C^1$ norm. 
Then
\[
  \liminf_{i\to\infty}\dimh(\Pi_{\mathcal F_i})_*\mu_i\geq \dimh(\Pi_{\mathcal F})_*\mu.\tag*{\qed}
\]
\end{cor}

Now we are going to show that weak convergence alone is not enough to get lower semicontinuity, unlike in the quasi-Bernoulli case with a common quasi-Bernoulli constant. The very same example will also show that the inequality is sharp in Theorem \ref{thm:lsc-joint}.

\begin{ex}Let $\mu$ be the Bernoulli measure on $\Sigma=\{0,1\}^{\N}$ with the probability vector $(1/2,1/2)$. Then $\mu([\mathbf u])=2^{-k}$ if $\mathbf u\in\Sigma_k$.

For each $i\in\N$, let $W_i$ be the word which we get when we concatenate all $2^i$ binary words of length $i$ in lexicographic order. E.g.
\[W_2=00011011.\]
Concatenate $W_i$ infinitely, to get
\[\omega^{(i)}=W_iW_iW_i\dots\]
Since $W_i$ has length $L_i=i2^i$, set
\[\mu_i=\frac{1}{L_i}\sum_{k=0}^{L_i-1}\delta_{\sigma^k\omega^{(i)}}.\]
In other words, we take whole shift orbit $\omega^{(i)}, \sigma\omega^{(i)}, \sigma^2\omega^{(i)},\dots,\sigma^{L_i-1}\omega^{(i)}$, and we distribute equal mass among the points.

Clearly $\mu_i$ is $\sigma$-invariant, because $\sigma$ permutes the points and they have equal mass.

Furthermore $\mu_i$ is ergodic, because if an invariant set contains one point from the orbit, then it contains all, i.e. it has either $0$ or $1$ mass.

Now we show that $\mu_i\Rightarrow\mu$. Fix $\mathbf u\in\Sigma_k$, and take $i>k$. Choose a starting position in $W_i$ uniformly. It means that we choose one of the $2^i$ blocks uniformly and then one of its $i$ positions uniformly. If the chosen position is among the first $i-k+1$ positions of a block, then the next $k$ digits lie within the block completely. Hence $2^{i-k}$ of the $2^i$ blocks have $\mathbf u$ at that position. Only the last $k-1$ positions can cross the boundary of a block. Therefore
\[\frac{i-k+1}{i}2^{-k}\leq\mu_i([\mathbf u])\leq\frac{(i-k+1)2^{i-k}+(k-1)2^i}{i2^i}=\frac{i-k+1}{i}2^{-k}+\frac{k-1}{i},\]
which gives that $\lim_i\mu_i([\mathbf u])=2^{-k}$.

Fix the following IFS on $[-1,2]$:
\[f_0(x)=\frac{x}{2},\qquad f_1(x)=\frac{x}{2}+\frac{1}{2}.\]
Evidently $\Pi_\cF(\omega)=\sum_{n=1}^\infty\omega_n2^{-n}$. As the support of $\mu_i$ consists of finitely many points, the support of $(\Pi_\cF)_*\mu_i$ is also finite, thus
\[\dimh(\Pi_\cF)_*\mu_i=0\quad (i\in\N).\]
But $(\Pi_\cF)_*\mu$ is the Lebesgue measure on $[0,1]$, because every dyadic interval of length $2^{-k}$ has mass $2^{-k}$. Therefore we get that weak convergence alone is not enough to get lower semicontinuity:
\[\liminf_i\dimh(\Pi_\cF)_*\mu_i=0<1=\dimh(\Pi_\cF)_*\mu.\]
Furthermore it can be readily seen that
\[h_{\mu_i}(\sigma)=0\,\,(i\in\N),\qquad h_{\mu}(\sigma)=\log 2,\qquad \chi_\mu(\cF)=\log 2.\]
If we plug in those values into \eqref{eq:hausdorff-dimension-map3}, we get 
\[0 = 1-\frac{\log 2-0}{\log 2}\]
showing that the inequality is sharp.
\qed
\end{ex}

Using Theorem \ref{thm:lsc-joint}, we can weaken both convergence conditions of our earlier joint lower semicontinuity result on Gibbs measures.

\begin{thm}For $0<\gamma_1<\gamma_2<1$ and $m\geq 2$, let $\mathcal F,\mathcal F_i\in\Theta_{\gamma_1,\gamma_2}(I)$ and $\mathcal F_i\to\mathcal F$ in $C^1$. Let $\varphi_i,\varphi$ be Hölder continuous potentials on $\Sigma$ (not necessarily in the same $C^\alpha$ Hölder class) such that $\varphi_i\to\varphi$ uniformly. Let $\mu_{\varphi_i},\mu_\varphi$ denote the (uniquely determined) associated shift-invariant Gibbs measures. Then
\begin{equation}
  \liminf_{i\to\infty}\dimh(\Pi_{\mathcal F_i})_*\mu_{\varphi_i}\geq\dimh(\Pi_{\mathcal F})_*\mu_\varphi.
\end{equation}
\end{thm}
\begin{proof}Write $\mu_i=\mu_{\varphi_i},\,\mu=\mu_\varphi$.

First we prove that $\mu_i\Rightarrow\mu$. Take a weak limit of $(\mu_i)$: $\mu_{i_j}\Rightarrow\nu$. A weak accumulation point exists, because the space of probability measures on $\Sigma$ is weakly compact. Since each $\mu_{i_j}$ is $\sigma$-invariant and $\sigma$ is continuous, $\nu$ is $\sigma$-invariant as well.

By the upper semicontinuity of the symbolic entropy (\cite[Theorem 8.2]{Walters82}),
\[h_\nu(\sigma)+\int\varphi\,d\nu\geq\limsup_j\left(h_{\mu_{i_j}}(\sigma)+\int\varphi_{i_j}\,d\mu_{i_j}\right)\]
\[=\limsup_jP(\varphi_{i_j})=\lim_jP(\varphi_{i_j})=P(\varphi).\]
Where we used that each \(\mu_{i_j}\) is an equilibrium state for its own potential \(\varphi_{i_j}\) (see \eqref{eq:vp-maximizer}), and $P(\varphi_{i_j})\to P(\varphi)$ holds because
\[|P(\varphi_{i_j})-P(\varphi)|\leq||\varphi_{i_j}-\varphi||_\infty.\]
By the variational principle (see \eqref{eq:variational-principle}) we get the reverse inequality
\[h_\nu(\sigma)+\int\varphi\,d\nu\leq P(\varphi).\]
Thus $\nu=\mu$ by the uniqueness of the shift-invariant Gibbs measure $\mu$. Since every weak accumulation point equals to $\mu$, we obtain $\mu_i\Rightarrow\mu$.

Now $\varphi_i\to\varphi$ uniformly and $\mu_i\Rightarrow\mu$ implies
\[h_{\mu_i}(\sigma)=P(\varphi_i)-\int\varphi_i\,d\mu_i\longrightarrow P(\varphi)-\int\varphi\,d\mu=h_\mu(\sigma).\]
Then by Theorem \ref{thm:lsc-joint} we get the claim.
\end{proof}

%---------------------------------------------------------------------------------------------------------------------------------------------------------Acknowledgment-----------------
\section*{Acknowledgments.}
I am grateful to professor Károly Simon, for his continuous support and encouragement throughout this research.

%---------------------------------------------------------------------------------------------------------------------------------------------------------thebibliography-----------------

{\footnotesize

%------------------------------------------------------------------------------------------------------------------------------------------------------------address---------------
\noindent HUN-REN–BME Stochastics Research Group, Budapest University of Technology and Economics, Műegyetem rkp. 3, H-1111 Budapest, Hungary

\noindent Budapest Semesters in Mathematics College International, Bethlen Gábor tér 2, Budapest, Hungary

\noindent email: alosonczi2@gmail.com

}
\end{document}